\documentclass[12pt]{amsart}
\usepackage{amsfonts}
\usepackage{pstricks,pstcol,pst-plot,pst-node,pst-tree,amssymb}
\usepackage{amsmath,amscd}
\usepackage{pstricks-add}
\usepackage{graphicx}

\newtheorem{thm}{\sc Theorem}[section]

\theoremstyle{definition}

\theoremstyle{definition}
\newtheorem{de}[thm]{\sc Definition}
\theoremstyle{definition}
\newtheorem{rem}[thm]{\sc Remark}
\theoremstyle{definition}

\numberwithin{equation}{section}
\input{tcilatex}

\begin{document}
\title[Spherical $CR$ manifolds]{Uniformization of spherical $CR$ manifolds}
\author{Jih-Hsin Cheng}
\address[Jih-Hsin Cheng]{Institute of Mathematics, Academia Sinica, Taipei,
11529 and National Center for Theoretical Sciences, Taipei Office, Taiwan,
R.O.C.}
\email{cheng@math.sinica.edu.tw}
\author{Hung-Lin~Chiu}
\address[Chiu]{Department of Mathematics, National Central University, Chung
Li, 32054, Taiwan, R.O.C.}
\email{hlchiu@math.ncu.edu.tw}
\author{Paul Yang}
\address[Yang]{Department of Mathematics, Princeton University, Princeton,
NJ 08544, U.S.A.}
\email{yang@Math.Princeton.EDU}

\begin{abstract}
Let $M$ be a closed (compact with no boundary) spherical $CR$ manifold of
dimension $2n+1$. Let $\widetilde{M}$ be the universal covering of $M.$ Let $%
\Phi $ denote a $CR$ developing map 
\begin{equation*}
\Phi :\widetilde{M}\rightarrow S^{2n+1}
\end{equation*}%
where $S^{2n+1}$ is the standard unit sphere in complex $n+1$-space $C^{n+1}$%
. Suppose that the $CR$ Yamabe invariant of $M$ is positive. Then we show
that $\Phi $ is injective for $n\geq 3$. In the case $n=2$, we also show
that $\Phi $ is injective under the condition: $s(M)<1$. It then follows
that $M$ is uniformizable.
\end{abstract}

\maketitle

\renewcommand{\subjclassname}{\textup{2000} Mathematics Subject 
Classification}

\section{Introduction and statement of the results}

Spherical $CR$ structures are modeled on the boundary of complex hyperbolic
space. There have been many studies in various aspects for this structure
(e.g., \cite{BS}, \cite{KT}, \cite{FG}, \cite{Gol}, \cite{CT}, \cite{Sch}).
In this paper, we study the uniformization problem. Let $S^{2n+1}$ denote
the standard unit sphere in complex $n+1$-space $C^{n+1}$. Let us start with
the $CR$ automorphism group $Aut_{CR}(S^{2n+1})$ of $S^{2n+1},$ which is the
group of complex fractional linear transformation $SU(n+1,1)/(\text{center})$%
. We have the following complex analogue of the Liouville theorem in
conformal geometry (\cite{CM}).

\bigskip

\textbf{Lemma 1.1.} Let $f$ be a $CR$ diffeomorphism from a connected open
set $U$ in $S^{2n+1}.$ If $f(U)\subset S^{2n+1}$, then $f$ is the
restriction to $U$ of a complex fractional linear transformation.

\bigskip

Let $M$ be a spherical $CR$ manifold of dimension $2n+1$. Let $\widetilde{M}$
be the universal covering of $M$. Using analytic continuation and Lemma 1.1,
we gets a $CR$ immersion $\Phi :\widetilde{M}\rightarrow S^{2n+1}.$ The map $%
\Phi $ is unique up to composition with elements of $Aut_{CR}(S^{2n+1})$
acting on $S^{2n+1}$. Such a map $\Phi $ is called a $CR$ developing map.

We will determine when $\Phi $ is injective. Let $\lambda (M)$ denote the $%
CR $ Yamabe invariant of $M$ (see (\ref{Ya}) in Section 2). In the case of $n
$ $=$ $2,$ we also need a condition on another $CR$ invariant $s(M)$ which
measures the integrability of a positive minimal Green's function $G_{p}$ on 
$\widetilde{M}$ (see Theorem 3.4 in Section 3):%
\begin{equation*}
s(M):=\inf {\{s:\int_{\widetilde{M}\setminus U_{p}}G_{p}^{s}dV_{\theta
}<\infty \}}
\end{equation*}%
\noindent where $U_{p}$ is a neighborhood of $p$.(see (\ref{s(M)}) in
Section 3). Observe that $s(M)$ $\leq $ $1$ in general (see Theorem 3.5 in
Section 3). We have the following result.

\bigskip

\textbf{Theorem A.} \textit{Let }$M$\textit{\ be a closed (compact with no
boundary) spherical }$CR$\textit{\ manifold of dimension }$2n+1$\textit{\
with\ }$\lambda (M)>0.$\textit{. Let }$\widetilde{M}$\textit{\ be the
universal covering of }$M.$\textit{\ Let }$\Phi $\textit{\ denote a }$CR$%
\textit{\ developing map }%
\begin{equation*}
\Phi :\widetilde{M}\rightarrow S^{2n+1}
\end{equation*}%
\noindent \textit{Then\ }$\Phi $ \textit{is injective for }$n\geq 3$\textit{%
. In case }$n=2$\textit{, }$\Phi $\textit{\ is injective if we further
assume }$s(M)<1$\textit{.}

\textit{\bigskip }

Theorem A implies that a closed spherical $CR$ manifold $M$ with $\lambda
(M)>0$ is uniformizable. Let $\pi _{1}(M)$ denote the fundamental group of $%
M.$ The $CR$ developing map $\Phi $ induces a group homomorphism: 
\begin{equation*}
\Phi _{\ast }:\pi _{1}(M)\rightarrow Aut_{CR}(S^{2n+1}).
\end{equation*}%
In case $\Phi $ is injective, the group homomorphism $\Phi _{\ast }$ is
injective. Note that $\pi _{1}(M)$ acts on $\widetilde{M}$ by deck
transformations. The following result follows from Theorem A.

\bigskip

\textbf{Corollary B.} \textit{Suppose that we are in the situation of
Theorem A. Then }$M$\textit{\ is }$CR$\textit{\ diffeomorphic to the
quotient }$\Omega /\Gamma $\textit{\ where }$\Omega =\Phi (\widetilde{M})$%
\textit{\ }$\subset $\textit{\ }$S^{2n+1}$\textit{\ and }$\Gamma =\Phi
_{\ast }(\pi _{1}(M))$ \textit{for }$n$ $\geq $ $2$\textit{. Moreover, }$%
\Gamma $\textit{\ is a discrete subgroup of }$Aut_{CR}(S^{2n+1})$\textit{\
and acts on }$\Omega $\textit{\ properly discontinuously.}

\textit{\bigskip }

The idea of the proof\ of Theorem A follows a similar line as for the
conformal case. Basically we will be dealing with the Green's functions of
the $CR$ invariant sublaplacian (see (\ref{ISL}) in Section 2) on different
spaces. In particular, the idea of comparing the pull-back $\bar{G}$ of the
Green's function on $S^{2n+1}$ with the (minimal positive) Green's function $%
G$ of $\widetilde{M}$ follows the work of Schoen and Yau (\cite{SY1} or \cite%
{SY2}). We reduce the injectivity problem to the estimate of the quotient $v$
$:=$ $\frac{G}{\bar{G}}.$ As expected, the $CR$ Bochner formula for $v$
contains an extra cross term which has no Riemannian analogue. Fortunately
we can manage this extra cross term by converting it into a term involving a
Paneitz-like operator $P$ (see (\ref{Pan}) in Section 2). The nonnegativity
of $P$ for $n$ $\geq $ $2$ (see line 5 in the proof of Proposition 3.2 in 
\cite{GL} or \cite{CC}) helps to simplify the estimates (see (\ref{24}) and (%
\ref{4}) in Section 4). We can finally prove $v$ $=$ $1,$ and hence $\Phi $
is injective under the condition mentioned in Theorem A.

There has been an unpublished paper (\cite{Li1}) about this uniformization
problem, circulating for years. The main difference between our paper and 
\cite{Li1} is the treatment of the $CR$ Bochner formula. We have realized
the important role of that Paneitz-like operator $P$ in the $CR$ setting of
the Bochner formula through the study of some other problems in recent years
(e.g., \cite{Chiu}, \cite{CC}). So we can clarify some estimates in \cite%
{Li1} and conclude a new result in the case of $n$ $=$ $2.$

Based on the uniformization of spherical $CR$ manifolds, in his another
unpublished paper (\cite{Li3}), using Bony's strong maximum principle, Z. Li
showed the nonnegativity of the CR mass (see Definition 5.1 in Section \ref%
{seccrmass}). We state his result as Corollary C:

\bigskip

\textbf{Corollary C.} \textit{\ Let} $M$ \textit{be a closed spherical }$%
\mathit{CR}$\textit{\ manifold with $\lambda (M)>0$.} \textit{Then, for $%
n\geq 3$, the\ }$CR$ \textit{mass\ }$A_{b}>0$ \textit{\ unless} $M$ \textit{%
\ is the standard sphere.} \textit{\ In case} $n=2,\ $ \textit{the same
result also holds if we assume futher } $s(M)<1$. \textit{\bigskip }

\begin{rem}
Notice that Z. Li's arguments are valid, provided that the $CR$ developing
map $\Phi $ is injective. We rewrite his proof in Section \ref{seccrmass}.
To solve the $CR$ Yamabe problem by producing a minimizer for the Yamabe (or
Tanaka-Webster) quotient, one needs to work out a test function estimate by
using the above positive mass theorem. This has been done in \cite{Li3} (see
also unpublished notes of Andrea Malchiodi). In dimension 3, we also have a
positive mass theorem under the condition $P$ $\geq $ $0$ through a
different approach (\cite{CMY}). Note that in dimension 3 (with $\lambda (M)$
$>$ $0)$, the condition $P$ $\geq $ $0$ is not automatic and is shown to be
almost equivalent to the embeddability of the underlying $CR$ structure (%
\cite{CCY1}, \cite{CCY2}).\newline
\end{rem}

\textbf{Acknowledgments.} The first author's research was supported in part
by NSC 98-2115-M-001-008-MY3, the second author's research was supported in
part by CIZE Foundation and in part by NSC 96-2115-M-008-017-MY3, and the
third author's research was supported in part by DMS-0758601.

\bigskip

\section{Bochner formulas and $CR$ invariant operators}

Let $M$ be a smooth (meaning $C^{\infty }$ throughout the paper$)$ ($2n+1)$%
-dimensional (paracompact) manifold. A contact structure or bundle $\xi $ on 
$M$ is a completely nonintegrable $2n$-dimensional distribution. A contact
form is a 1-form annihilating $\xi $. Let $(M,\xi )$ be a contact ($2n+1)$%
-dimensional manifold with an oriented contact structure $\xi $. There
always exists a global oriented contact form $\theta $, obtained by patching
together local ones with a partition of unity. The Reeb vector field of $%
\theta $ is the unique vector field $T$ such that $\theta (T)=1$ and $%
\mathcal{L}_{T}\theta =0$ or $d\theta (T,{\cdot })=0$. A $CR$-structure
compatible with $\xi $ is a smooth endomorphism $J:{\xi }$ ${\rightarrow }$ $%
{\xi }$ such that $J^{2}=-Identity$. Let $T_{1,0}$ $\subset $ $\xi \otimes C$
denote the $n$-dimensional complex subbundle of $TM\otimes C,$ consisting of
eigenvectors of $J$ with eigenvalue $i.$ We will assume throughout that the $%
CR$ structure $J$ is integrable, that is, $T_{1,0}$ satisfies the condition $%
[T_{1,0},$ $T_{1,0}]$ $\subset $ $T_{1,0}.$ A pseudohermitian structure
compatible with an oriented contact structure $\xi $ is a $CR$-structure $J$
compatible with $\xi $ together with a global contact form $\theta $. On $%
\xi ,$ we define the Levi form $L_{\theta }$ $:=$ $\frac{1}{2}d\theta (\cdot
,J\cdot ).$ If $L_{\theta }$ is definite (independent of the choice of
contact form), $M$ (or $(M,\xi ,J)$) is said to be strictly pseudoconvex. We
call $\theta $ positive if $L_{\theta }$ is positive definite (often called
Levi metric in this case). We will always assume that $M$ is strictly
pseudoconvex and $\theta $ is positive.

Given a pseudohermitian structure $(J,\theta )$ (with $J$ integrable and $%
\theta $ positive), we can choose complex vector field $Z_{\alpha }$, $%
\alpha $ $=$ $1,$ $2,...,$ $n,$ eigenvectors of $J$ with eigenvalue $i$, and
complex 1-form ${\theta }^{\alpha },$ $\alpha $ $=$ $1,$ $2,...,$ $n,$ such
that $\{\theta ,{\theta ^{\alpha }},{\theta ^{\bar{\alpha}}}\}$ is dual to $%
\{T,Z_{\alpha },Z_{\bar{\alpha}}\}$ (${\theta ^{\bar{\alpha}}}={\bar{({%
\theta ^{\alpha }})}}$,$Z_{\bar{\alpha}}={\bar{({Z_{\alpha }})}}$). It
follows that 
\begin{equation*}
d\theta =ih_{\alpha \bar{\beta}}{\theta ^{\alpha }}{\wedge }{\theta }^{\bar{%
\beta}}
\end{equation*}%
\noindent (summation convention throughout) for some hermitian matrix of
functions $(h_{\alpha {\bar{\beta}}})$, which is positive definite since $M$
is strictly pseudoconvex and $\theta $ is positive.

The pseudohermitian connection of $(J,\theta )$ is the connection $\nabla
^{p.h.}$ on $TM{\otimes }C$ (and extended to tensors) given by

\begin{equation*}
{\nabla }^{p.h.}Z_{\alpha }={\omega _{\alpha }}^{\beta }{\otimes }Z_{\beta },%
{\nabla }^{p.h.}Z_{\bar{\alpha}}={\omega _{\bar{\alpha}}}^{\bar{\beta}}{%
\otimes }Z_{\bar{\beta}},{\nabla }^{p.h.}T=0
\end{equation*}

\noindent in which the 1-forms ${\omega _{\alpha }}^{\beta }$ are uniquely
determined by the following equations with a normalization condition (\cite%
{We}, \cite{Ta}, \cite{Lee}):%
\begin{eqnarray}
d{\theta ^{\beta }} &=&{\theta ^{\alpha }}{\wedge }{\omega _{\alpha }}%
^{\beta }+{A^{\beta }}_{\bar{\alpha}}\theta {\wedge }{\theta ^{\bar{\alpha}},%
}  \label{2.1} \\
dh_{\alpha \bar{\gamma}} &=&{\omega _{\alpha }}^{\beta }h_{\beta \bar{\gamma}%
}+h_{\alpha \bar{\beta}}{\omega _{\bar{\gamma}}}^{\bar{\beta}}.  \notag
\end{eqnarray}

\noindent The coefficient ${A^{\beta }}_{\bar{\alpha}}$ in ($\ref{2.1})$ is
called the (pseudohermitian) torsion. As usual we use the matrix $h_{\alpha 
\bar{\beta}}$ to raise and lower indices, e.g., ${A}_{\alpha \gamma }$ $=$ $%
h_{\alpha \bar{\beta}}{A^{\bar{\beta}}}_{\gamma }$ where ${A^{\bar{\beta}}}%
_{\gamma }$ is the complex conjugate of ${A^{\beta }}_{\bar{\gamma}}.$ We
define covariant differentiation with respect to the connection $\nabla
^{p.h.}.$ For a real $C^{\infty }$ smooth function $\varphi ,$ we have $%
\varphi _{0}$ $:=$ $T\varphi ,$ $\varphi _{\alpha }$ $:=$ $Z_{\alpha
}\varphi $, $\varphi _{\alpha \beta }$ $:=$ $Z_{\beta }Z_{\alpha }\varphi -{%
\omega _{\alpha }}^{\gamma }(Z_{\beta })Z_{\gamma }\varphi ,$ etc.. For the
subgradient $\nabla _{b}\varphi $, the sublaplacian ${\Delta }_{b}\varphi ,$
and the subhessian ($\nabla ^{H}$)$^{2}\varphi $, we have the following
formulas:%
\begin{eqnarray*}
\nabla _{b}\varphi &=&\varphi ^{\alpha }Z_{\alpha }+\varphi ^{\bar{\alpha}%
}Z_{\bar{\alpha}} \\
{\Delta }_{b}\varphi &=&-(\varphi {_{\alpha }}^{\alpha }+\varphi {_{\bar{%
\alpha}}}^{\bar{\alpha}}) \\
(\nabla ^{H})^{2}\varphi &=&\varphi {_{\alpha }}^{\beta }\theta ^{\alpha
}\otimes Z_{\beta }+\varphi {_{\bar{\alpha}}}^{\beta }\theta ^{\bar{\alpha}%
}\otimes Z_{\beta } \\
&&+\varphi {_{\alpha }}^{\bar{\beta}}\theta ^{\alpha }\otimes Z_{\bar{\beta}%
}+\varphi {_{\bar{\alpha}}}^{\bar{\beta}}\theta ^{\bar{\alpha}}\otimes Z_{%
\bar{\beta}}.
\end{eqnarray*}%
Differentiating ${\omega _{\beta }}^{\alpha }$ gives

\begin{eqnarray*}
&&d{\omega _{\beta }}^{\alpha }{-\omega _{\beta }}^{\gamma }\wedge {\omega
_{\gamma }}^{\alpha } \\
&=&R{_{\beta }}^{\alpha }{}_{\rho \bar{\sigma}}\theta ^{\rho }\wedge \theta
^{\bar{\sigma}}+i{A^{\alpha }}_{\bar{\gamma}}\theta _{\beta }\wedge \theta ^{%
\bar{\gamma}}-iA_{\beta \gamma }\theta ^{\gamma }\wedge \theta ^{\alpha }%
\text{ mod }\theta
\end{eqnarray*}

\noindent where $R{_{\beta }}^{\alpha }{}_{\rho \bar{\sigma}}$is the
Tanaka-Webster curvature. Write $R_{\alpha \bar{\beta}}$ $:=$ $R{_{\gamma }}%
^{\gamma }{}_{\alpha \bar{\beta}}$ and $R$ $:=$ $R{_{\alpha }}^{\alpha }.$
For $X$ $=$ $X^{\alpha }Z_{\alpha },$ $Y$ $=$ $Y^{\beta }Z_{\beta }$ $\in $ $%
T_{1,0},$ we define 
\begin{eqnarray*}
Ric(X,Y) &=&R_{\alpha \bar{\beta}}X^{\alpha }Y^{\bar{\beta}} \\
Tor(X,Y) &=&2\func{Re}(iA_{\bar{\alpha}\bar{\beta}}X^{\bar{\alpha}}Y^{\bar{%
\beta}}).
\end{eqnarray*}

We recall the pointwise Bochner formula (\cite{Gr}):%
\begin{eqnarray}
\frac{1}{2}{\Delta }_{b}|\nabla _{b}\varphi |^{2} &=&-|(\nabla
^{H})^{2}\varphi |^{2}+<\nabla _{b}\varphi ,\nabla _{b}{\Delta }_{b}\varphi >
\label{BF} \\
&&-2Ric((\nabla _{b}\varphi )_{C},(\nabla _{b}\varphi )_{C})  \notag \\
&&+(n-2)Tor((\nabla _{b}\varphi )_{C},(\nabla _{b}\varphi )_{C})  \notag \\
-2 &<&J\nabla _{b}\varphi ,\nabla _{b}\varphi _{0}>  \notag
\end{eqnarray}

\noindent for a real smooth function $\varphi ,$ where the length 
\TEXTsymbol{\vert} $\cdot $ $|$ and the inner product $<\cdot ,\cdot >$ are
with respect to the Levi metric $L_{\theta }$ and $(\nabla _{b}\varphi )_{C}$%
$:=$ $\varphi ^{\alpha }Z_{\alpha }.$

We define a Paneitz-like operator $P$ by%
\begin{equation}
P\varphi :=4(\varphi {_{\bar{\alpha}}}^{\bar{\alpha}}{_{\beta }}+inA_{\beta
\alpha }\varphi ^{\alpha })^{\beta }+conjugate.  \label{Pan}
\end{equation}

\noindent Let $P_{\beta }\varphi $ $:=$ $\varphi {_{\bar{\alpha}}}^{\bar{%
\alpha}}{_{\beta }}+inA_{\beta \alpha }\varphi ^{\alpha }.$ For $n$ $=$ $1$,
the $CR$ pluriharmonic functions are characterized by $P_{1}\varphi $ $=$ $0$
(\cite{Lee2}) while, for $n$ $\geq $ $2,$ they are characterized by $%
P\varphi $ $=$ $0.$ (see \cite{GL} in which $P$ is also identified with the
compatibility operator for solving a certain degenerate Laplace equation in
the case of $n$ $=$ $1$). On the other hand, this operator $P$ is a $CR$
analogue of the Paneitz operator in conformal geometry (see \cite{Hir} for
the relation to a $CR$ analogue of the $Q$-curvature and the $\log $-term
coefficient in the Szeg\"{o} kernel expansion). On a closed pseudohermitian (%
$2n+1)$-dimensional manifold $(M,$ $J,$ $\theta ),$ we call $P$ nonnegative
if there holds%
\begin{equation}
\int_{M}\varphi (P\varphi )dV_{\theta }\geq 0  \label{NG}
\end{equation}

\noindent for all real smooth functions $\varphi ,$ in which the volume form 
$dV_{\theta }$ $=$ $\theta \wedge (d\theta )^{n}.$ We need the integrated
Bochner formula:%
\begin{eqnarray}
\int_{M}\varphi _{0}^{2}dV_{\theta } &=&\frac{1}{n^{2}}\int_{M}({\Delta }%
_{b}\varphi )^{2}dV_{\theta }  \label{IBF} \\
&&+\frac{2}{n}\int_{M}Tor((\nabla _{b}\varphi )_{C},(\nabla _{b}\varphi
)_{C})dV_{\theta }  \notag \\
&&-\frac{1}{2n^{2}}\int_{M}\varphi (P\varphi )dV_{\theta }  \notag
\end{eqnarray}

\noindent (see Corollary 2.4 in \cite{CC}). By Theorem 3.2 in \cite{CC}, we
learn that $P$ is nonnegative for $n$ $\geq $ $2.$ It follows from (\ref{IBF}%
) and (\ref{NG}) that%
\begin{eqnarray}
\int_{M}\varphi _{0}^{2}dV_{\theta } &\leq &\frac{1}{n^{2}}\int_{M}({\Delta }%
_{b}\varphi )^{2}dV_{\theta }  \label{IBF2} \\
&&+2\kappa \int_{M}|\nabla _{b}\varphi |^{2}dV_{\theta }  \notag
\end{eqnarray}

\noindent where $\kappa $ $:=$ $\max_{q\in M}(\sum_{\alpha ,\beta
}(A_{\alpha \beta }A^{\alpha \beta })(q))^{1/2}$ (note that $\sum_{\alpha
,\beta }A_{\alpha \beta }A^{\alpha \beta }$ is independent of the choice of
frames and is a nonnegative real function on $M$).

Let $b_{n}$ $:=$ $2+\frac{2}{n}.$ We define the $CR$ invariant sublaplacian $%
D_{\theta }$ by%
\begin{equation}
D_{\theta }=b_{n}{\Delta }_{b}+R  \label{ISL}
\end{equation}

\noindent where $R$ denotes the Tanaka-Webster scalar curvature (with
respect to $\theta $ while fixing $J)$. Suppose that $\tilde{\theta}$ $=$ $%
u^{\frac{2}{n}}\theta $ for a positive $C^{\infty }$ smooth function $u.$
Then for any real smooth function $\varphi ,$ there holds%
\begin{equation}
D_{\theta }(u\varphi )=u^{1+\frac{2}{n}}D_{\tilde{\theta}}(\varphi ).
\label{TL}
\end{equation}

\noindent Letting $\varphi $ $\equiv $ $1$ in (\ref{TL}) gives the
transformation law for $R:$%
\begin{equation*}
\tilde{R}=u^{-1-\frac{2}{n}}D_{\theta }(u)
\end{equation*}

\noindent where $\tilde{R}$ denotes the Tanaka-Webster scalar curvature with
respect to $\tilde{\theta}.$ The Yamabe problem on a $CR$ manifold is to
find $u$ (or $\tilde{\theta})$ such that $\tilde{R}$ is a given constant.
This is the Euler-Lagrange equation of the following energy functional:%
\begin{equation*}
E_{\theta }(u):=\int_{M}(b_{n}|\nabla _{b}u|^{2}+Ru^{2})dV_{\theta }.
\end{equation*}

\noindent for positive smooth functions $u$ such that 
\begin{equation}
\int_{M}|u|^{b_{n}}dV_{\theta }=1.  \label{C}
\end{equation}

\noindent The $CR$ Yamabe invariant $\lambda (M)$ has the following
expression:%
\begin{equation}
\lambda (M)=\inf_{u\in \Xi }E_{\theta }(u)  \label{Ya}
\end{equation}

\noindent where $\Xi $ is the space of positive smooth (with compact support
if $M$ is noncompact) functions $u$ satisfying (\ref{C}). For $M$ closed, it
is known that $\lambda (M)$ $>$ $0$ is equivalent to the existence of a
contact form $\bar{\theta}$ with respect to which $\bar{R}$ $>$ $0.$

\bigskip

Let $(M,$ $J,$ $\theta )$ be a closed pseudohermitian manifold with $R>0.$
Let $\Gamma ^{\beta }(M)$ denote the nonisotropic H\"{o}lder space of
exponent $\beta $ (page 181 in \cite{JL} or \cite{FS} for the local
description modelled on the Heisenberg group). Following a standard argument
in \cite{A} for the elliptic case, we obtain

\bigskip

\textbf{Proposition 2.1.} \textit{Let }$(M,$\textit{\ }$J,$\textit{\ }$%
\theta )$\textit{\ be a closed pseudohermitian manifold with }$R>0.$\textit{%
\ Then for any }$f\in \Gamma ^{\beta }(M),\ \beta $\textit{\ a noninteger }$%
> $\textit{\ }$0$\textit{, there exists a unique }$u\in \Gamma ^{\beta
+2}(M) $\textit{\ such that }%
\begin{equation*}
D_{\theta }(u)=f.
\end{equation*}

\bigskip

Using Proposition 2.1 and a similar construction in \cite{A}, we have that
for any $p\in M$, there is a unique Green's function $G_{p}$ for $D_{\theta
} $ with pole at $p$.

\bigskip

\section{The Green's function of the universal covering}

Let $S^{2n+1}$ denote the unit sphere in $C^{n+1}.$ Let $\hat{\xi}$ $:=$ $%
TS^{2n+1}$ $\cap $ $J_{C^{n+1}}$ $(TS^{2n+1})$ be the standard contact
bundle over $S^{2n+1},$ where $J_{C^{n+1}}$ denotes the almost complex
structure of $C^{n+1}.$ Let $\hat{J}$ be the restriction of $J_{C^{n+1}}$ to 
$\hat{\xi}.$ We call a $CR$ manifold $(M,$ $J)$ (or $(M,$ $\xi ,$ $J),$ resp.%
$)$ spherical if it is locally $CR$ equivalent to ($S^{2n+1},$ $\hat{J})$
(or $(S^{2n+1},$ $\hat{\xi},$ $\hat{J}),$ resp.) (cf. e.g., \cite{BS}). Let $%
(M,$ $J,$ $\theta )$ be a closed pseudohermitian manifold of dimension $2n+1$
with $(M,$ $J)$ spherical and $R$ $>$ $0.$ Let $\widetilde{M}$ be the
universal covering of $M$ with the $CR$ structure $\pi ^{\ast }J$ and the
contact form $\pi ^{\ast }\theta $, where 
\begin{equation*}
\pi :\widetilde{M}\rightarrow M
\end{equation*}%
is the canonical projection. It follows that $\widetilde{M}$ has no
boundary. If $\widetilde{M}$ is compact, then $(\widetilde{M},$ $\pi ^{\ast
}J)$ must be $CR$ equivalent to ($S^{2n+1},$ $\hat{J})$ since it is simply
conected and spherical. We will assume that $\widetilde{M}$ is noncompact
(or $\pi _{1}(M)$ is an infinite group) throughout this section. We will
still use $\theta $ to mean $\pi ^{\ast }\theta $. Our goal in this section
is to study the existence of the Green's function for $D_{\theta }$ on $%
\widetilde{M}$ and its decay property at the geometric boundary of $%
\widetilde{M}$.

Let $\Omega $ be a relatively compact smooth domain in $\widetilde{M}$ and $%
p $ $\in $ $\Omega $. We would like to construct the Dirichlet Green's
function $G_{p}^{\Omega }$ for the domain $\Omega $, that is, to prove the
following

\bigskip

\textbf{Theorem 3.1.} \textit{There exists a unique }$G_{p}^{\Omega }\in
C^{\infty }(\Omega \setminus \{p\})\cap C(\overline{\Omega }\setminus \{p\})$%
\textit{\ such that}%
\begin{equation}
\begin{array}{rcl}
D_{\theta }(G_{p}^{\Omega }) & = & \delta _{p}\ \text{in}\ \Omega \\ 
G_{p}^{\Omega }|_{\partial \Omega } & = & 0%
\end{array}
\label{3.1}
\end{equation}

\textit{\bigskip }

Once $G_{p}^{\Omega }$ is constructed, the symmetry property $G_{p}^{\Omega
}(q)$ $=$ $G_{q}^{\Omega }(p)$is due to the fact that $D_{\theta }$ is
self-adjoint and from the integration by parts argument as in the elliptic
case. The positivity of $G_{p}^{\Omega }$ is due to the fact that the
leading order operator of $D_{\theta }$ is nonnegative and the
Tanaka-Webster scalar curvature $R$ is positive.

As in the elliptic case, the existence of the Dirichlet Green's function is
equivalent to the solvability of nonhomogeneous Dirichlet problem with zero
boundary value. So solving (\ref{3.1}) is reduced to solving the following
Dirichlet problem..

\bigskip

\textbf{Theorem 3.2.} \textit{Let }$\Omega $\textit{\ be a relatively
compact smooth domain in }$\widetilde{M}$\textit{\ and }$f\in \Gamma ^{\beta
}(\overline{\Omega })$\textit{. Then there is a unique }$u\in \Gamma ^{\beta
+2}(\Omega )\cap C(\overline{\Omega })$\textit{\ such that}%
\begin{equation}
\begin{array}{rcl}
D_{\theta }(u) & = & f\ \text{in}\ \Omega \\ 
u|_{\partial \Omega } & = & 0%
\end{array}
\label{3.2}
\end{equation}

\textit{\bigskip }

The uniquness of the solution to (\ref{3.2}) follows from the following
lemma.

\bigskip

\textbf{Lemma 3.3.}\textit{\ Let }$\Omega $\textit{\ be a relatively compact
smooth domain in }$\widetilde{M}$\textit{\ such that for }$u,v\in
C^{2}(\Omega )\cap C(\overline{\Omega }),$\textit{\ there holds}%
\begin{eqnarray*}
D_{\theta }(u) &\leq &D_{\theta }(v)\text{ in}\ \Omega \\
u|_{\partial \Omega } &=&v|_{\partial \Omega }.
\end{eqnarray*}%
\noindent \textit{Then }$u$\textit{\ }$< $\textit{\ }$v$\textit{\ in }$%
\Omega $\textit{\ unless\ }$u=v$ \textit{\ in} $\Omega$.

\bigskip

To prove Lemma 3.3, we observe that the leading order part of $D_{\theta }$
is a subelliptic operator of Hormander's type (sum of square vector fields).
Then one can apply Bony's arguments without essential change by the local
nature of this lemma. To prove Theorem 3.2, we use Perron's construction
which relies heavily on the maximum principle and the solvability of the
Dirichlet problem for the balls. We remark that the key step of proving
Theorem 3.2 is to localize the problem. In this respect, the $CR$ invariance
enables us to reduce the problem on $(\widetilde{M},\theta )$ ($J$ omitted)
for $D_{\theta }$ to the problem on the Heisenberg group $(H^{n},\Theta )$
for $D_{\Theta }$ (\cite{JL}).

We will often call a Heisenberg ball simply a ball in this section. Let $B$
denote a small ball in $\widetilde{M}$, identified with a Heisenberg ball in 
$H^{n}.$ There is a positive function $\phi \in C^{\infty }(\widetilde{M})$
such that $\phi ^{\frac{2}{n}}\theta =\Theta $ in $B$. By the known results
on $H^{n}$ and the transformation law for $D_{\theta }$, we conclude that
for each $f\in \Gamma ^{\beta }(\bar{B})$ and $g\in C(\partial {B})$, there
is a unique $u\in \Gamma ^{\beta +2}(B)\cap C(\bar{B})$ such that%
\begin{eqnarray*}
D_{\theta }(u) &=&f\ \text{in}\ B \\
u|_{\partial B} &=&g
\end{eqnarray*}

\noindent In the Dirichlet problem for $D_{\theta }$ on a smooth domain in $%
\widetilde{M}$, the question of continuity up to the boundary is a purely
local issue. So we will deal with it in the same localizing spirit as above.

We begin the Perron process by generalizing the notion of subsolution in
classical elliptic theory to the operator $D_{\theta }$ on a smooth domain $%
\Omega $ in $\widetilde{M}$. Note that $D_{\theta }$ has a nonnegative
leading order part.

\bigskip

\textbf{Definition.} A continuous function $u$ in $\Omega $ is called a
subsolution to the equation 
\begin{equation*}
D_{\theta }(v)=f
\end{equation*}%
\noindent where $f\in \Gamma ^{\beta }(\bar{\Omega}),\beta $ a noninteger $>$
$0$, if for every ball $B\subset \subset \Omega $ and $v$ such that 
\begin{equation*}
D_{\theta }(v)=f\ \text{in}\ B,\ u\leq v\ \text{on}\ \partial {B},
\end{equation*}%
\noindent then we have that $u\leq v$ in $B$.

Analogously we can define the notion of supersolution as well. These notions
are completely in parallel to the notions of continuous subharmonic and
superharmonic functions in classical elliptic theory. The significance of
these notions are ensured by the Bony's maximum principle. They also have
the following useful properties of sub and supersolutions:

(1) If $u\in C^{2}(\Omega )$, then $u$ is a subsolution (supersolution,
resp.) if and only if that $D_{\theta }(u)\leq f$ $(D_{\theta }(u)\geq f,$
resp.$).$

(2) If $u_{1},\cdots ,u_{m}$ are subsolutions (supersolutions, resp.) in $%
\Omega $, then $\max \{u_{j}:1\leq j\leq m\}$ $(\min \{u_{j}:1\leq j\leq
m\}, $ resp.$)$ is also a subsolution (supersolution, resp.) in $\Omega $.

(3) Suppose that $B\subset \subset \Omega $ and $u_{1}$ satisfies 
\begin{eqnarray*}
D_{\theta }(u_{1}) &=&f\text{ in}\ B \\
u_{1} &=&u_{2}\text{ in}\ \partial B
\end{eqnarray*}%
\noindent where $u_{2}$ is a subsolution in $\Omega $. Then%
\begin{equation*}
u=\Big\{%
\begin{array}{ll}
u_{1} & \ \text{in}\ B \\ 
u_{2} & \ \text{in}\ \Omega \backslash B%
\end{array}%
\end{equation*}%
\noindent is also a subsolution in $\Omega $.

\bigskip

\proof
\textbf{(of Theorem 3.2)} We will carry out the proof in the spirit of
standard arguments in the elliptic theory. It consists of two main steps:

\bigskip

\textbf{Step $1$ : Construction of the Perron solution.} Consider the
folllowing set of subsolutions 
\begin{equation*}
S=\{v:v\ \text{is a subsolution for}\ D_{\theta }(w)=f,\text{ }v|_{\partial
\Omega }\leq 0\}
\end{equation*}%
Note that the Tanaka-Webster scalar curvature on $\widetilde{M}$ has a
positive lower bound, say $R_{0}$. We observe that $\frac{-\sup {|f|}}{R_{0}}%
\in S$ and the constant $\frac{\sup {|f|}}{R_{0}}$ is a supersolution.
Therefore $u(x)=\sup_{v\in S}{v(x)}$ is well defined.

We would like to show that $u$ is the Perron solution, i.e. ,$D_{\theta}u=f$
in $\Omega$. Let $p\in \Omega$ be an arbitrary fixed point. By the
definition of $u$, there exists a sequence $v_{m}\in S$ such that $%
v_{m}(p)\rightarrow u(p)$. By replacing $v_{m}$ with $\max\{v_{1}, \cdots,
v_{m}\}$, we may assume that the sequence is monotone.

Now choose a ball $B\subset \subset \Omega $ of $p$ such that the geometry
of $B$ can be flattened after a conformal change of a contact form. Let $%
w_{m}$ be the unique solution satisfying%
\begin{equation*}
\begin{array}{rcll}
D_{\theta }(w_{m}) & = & f & \text{in}\ B \\ 
w_{m} & = & v_{m} & \text{in}\ \partial B.%
\end{array}%
\end{equation*}%
\noindent It follows from earlier discussion that%
\begin{equation*}
W_{m}=\Big\{%
\begin{array}{ll}
w_{m} & \ \text{in}\ B \\ 
v_{m} & \ \text{in}\ \Omega \backslash B%
\end{array}%
\end{equation*}%
\noindent is a subsolution in $\Omega $. Hence, $u(p)\geq W_{m}(p)\geq
v_{m}(p)\rightarrow u(p)$. Because $W_{m}$ is monotone increasing in $B$,
the limit $W=\lim_{m\rightarrow \infty }W_{m}$ exists in $B$. As in the
elliptic theory, the subelliptic apriori estimates imply that the sequence $%
W_{m}$ contains a subsequence converging in $B$, and hence $W$ is a solution
in $B$ and $W(p)=u(p)$.

We claim that $W=u$ in $B$ to complete step $1$. The arguments are standard:
let $q\in B$, there is a monotone increasing sequence $g_{m}\in S$ such that 
$g_{q}\rightarrow u(q)$. Let $h_{m}$ solve the equation%
\begin{equation*}
\begin{array}{rcll}
D_{\theta }(w) & = & f & \text{in}\ B \\ 
w & = & \bar{g}_{m} & \text{in}\ \partial B%
\end{array}%
\end{equation*}

\noindent where in $B,\ \bar{g}_{m}=\max \{g_{m},v_{m}\}$. Therefore the
sequence is also monotone increasing and $v_{m}\leq h_{m}$ in $B$. As
before, $h=\lim_{m\rightarrow \infty }h_{m}$ is a solution in $B$ and $%
h(q)=u(q)$. Since $u(p)\geq h(p)\geq W(p)=u(p)$, Bony's strong maximum
principle implies that $W=h$ in $B$. Therefore, $u$ is indeed an interior
solution.

\bigskip

\textbf{Step $2$ : Continuity up to the boundary.} Fix $p\in \partial \Omega 
$, we choose a small ball $B$ of $p$ such that the boundary of $B\cap \Omega 
$ is smooth and the geometry of $B$ can be flattened by choosing a conformal
contact form. Notice the following two facts:

(1) There exists $u_{f}$ which solves the Dirichlet problem:%
\begin{equation*}
\begin{array}{rcll}
D_{\theta }(u_{f}) & = & f & \text{in}\ B\cap \Omega \\ 
u_{f} & = & 0 & \text{in}\ \partial (B\cap \Omega ).%
\end{array}%
\end{equation*}

(2) The existence of a local barrier at $p$, that is, for a small ball of $p$%
, there is a function $w\in C(\overline{B}\cap \overline{\Omega })\cap
C^{2}(B\cap \Omega )$ such that 
\begin{equation*}
D_{\theta }(w)=0\ \text{in}\ B\cap \Omega ,w>0\ \text{in}\ \overline{B}\cap 
\overline{\Omega }\backslash \{p\}\ \text{and}\ w(p)=0.
\end{equation*}

We will use this local barrier to construct a global barrier and show that
the Perron solution $u$ obtained in step 1 is continuous up to the boundary.
For any $\varepsilon $, there exists $\delta $ and $K$ such that 
\begin{equation*}
|u_{f}(x)|\leq \varepsilon \ \text{if}\ |x-p|\leq \delta
\end{equation*}%
\noindent and 
\begin{equation*}
Kw(x)\geq \sup {u_{f}}\ \text{if}\ |x-p|\geq \delta .
\end{equation*}%
\noindent We will consider $w_{1}=-Kw+u_{f}-\varepsilon $ in $B\cap \Omega $%
. Then we have immediately that 
\begin{equation*}
D_{\theta }w_{1}=f\ \text{and}\ w_{1}<0\ \text{in}\ B\cap \Omega .
\end{equation*}%
Consider the following number 
\begin{equation*}
M_{K}=\sup_{\partial (B\cap \Omega )\backslash \partial \Omega
}w_{1}=\sup_{\partial (B\cap \Omega )\backslash \partial \Omega
}(-Kw-\varepsilon )=-K(\min_{\partial (B\cap \Omega )\backslash \partial
\Omega }w)-\varepsilon ,
\end{equation*}%
\noindent and we let%
\begin{equation*}
W_{1}=\Big\{%
\begin{array}{ll}
\max (w_{1},M_{K}) & \ \text{in}\ \overline{B}\cap \overline{\Omega } \\ 
M_{K} & \ \text{in}\ \overline{\Omega }\backslash \overline{B}.%
\end{array}%
\end{equation*}
\noindent Note that $W_{1}|_{\partial (B\cap \Omega )\backslash \partial
\Omega }=M_{K}$, so $W_{1}\in C(\Omega )$. It is easy to check directly that 
$W_{1}$ is a subsolution (when $K$ is large enough). To check the boundary
behavior of $W_{1}$, we note that (when $q$ is very close to $p$): 
\begin{equation*}
W_{1}(q)=\max (W_{1}(q),M_{K})=-Kw(q)+u_{f}(q)-\varepsilon \rightarrow
-\varepsilon
\end{equation*}%
as $q\rightarrow p$. So it is zero at $p$ and negative everywhere else on $%
\partial \Omega $. Therefore, $W_{1}\in S$ and $W_{1}$ is a global boundary
barrier at $p$. It follows that 
\begin{equation*}
\lim \inf_{q\rightarrow p}u(q)\geq \lim_{q\rightarrow
p}w_{1}(q)=-\varepsilon .
\end{equation*}%
Since $\varepsilon $ is arbitrary, we have the continuity of $u$ up to the
boundary.

\endproof%

\bigskip

We are now ready to construct a Green's function $G_{p}$ for $D_{\theta }$
with pole at $p$ $\in $ $\widetilde{M}$. We would also like to discuss its
decay properties at the infinity.

Recall that $\Phi :\widetilde{M}\rightarrow S^{2n+1}$ denotes the $CR$
developing map. Let $H_{y}$ be the Green's function for the $CR$ invariant
sublaplacian $D_{0}$ of $(S^{2n+1},\theta _{S^{2n+1}})$ with the pole $%
y=\Phi (p)$, where $\theta _{S^{2n+1}}$ is the standard contact form on $%
S^{2n+1}$. Since $\Phi $ is a $CR$ immersion, we can write 
\begin{equation*}
\Phi ^{\ast }(\theta _{S^{2n+1}})=|\Phi ^{^{\prime }}|^{2}\theta
\end{equation*}%
\noindent where $\theta $ is the contact form for $\widetilde{M}$ and $|\Phi
^{^{\prime }}|$ is a positive $C^{\infty }$ smooth function on $\widetilde{M}
$. By the transformation law (\ref{TL}) of the $CR$ invariant sublaplacian,
we immediately obtain the following formula 
\begin{equation}
D_{\theta }(|\Phi ^{^{\prime }}|^{n}H_{y}\circ \Phi )=\sum_{\bar{p}\in \Phi
^{-1}(y)}|\Phi ^{^{\prime }}(\bar{p})|^{n+2}\delta _{\bar{p}}.
\label{3.11.1}
\end{equation}%
\noindent Let us consider the following function (with poles in $\Phi
^{-1}(y)$): 
\begin{equation}
\overline{G}:=|\Phi ^{^{\prime }}(p)|^{-(n+2)}|\Phi ^{^{\prime
}}|^{n}H_{y}\circ \Phi .  \label{3.11.2}
\end{equation}%
\noindent Since $\Phi $ is a $CR$ immersion, $\overline{G}$ is positive, $%
C^{\infty }$ smooth, and $D_{\theta }\overline{G}=0$ on $\widetilde{M}%
\backslash \Phi ^{-1}(y).$ Also, $\overline{G}$ has exactly the same
asymptotic behavior at each of $\Phi ^{-1}(y)$. We call $\overline{G}$ the
normalized pullback of $H_{y}$, which will be taken as a singular barrier in
the construction of a global Green's function on $\widetilde{M}$ through a
limit procedure of Dirichlet Green's functions.

Let $\{\Omega _{k}\subset \Omega _{k+1}:k=1,\cdots ,\}$ be a relatively
compact, $C^{\infty }$ smooth exhaustion of the universal covering $%
\widetilde{M}.$ Take $p$ $\in $ $\Omega _{1}$. Note that 
\begin{equation*}
D_{\theta }(\overline{G}-G_{p}^{\Omega _{k}})\geq 0,
\end{equation*}%
\noindent $\overline{G}$ is positive, and $G_{p}^{\Omega _{k}}|_{\partial
\Omega _{k}}=0$. By the maximum principle of Bony, we see that 
\begin{equation}
G_{p}^{\Omega _{k}}<\overline{G}.  \label{barrier}
\end{equation}%
\noindent away from $\Phi ^{-1}(y).$ Near the point $p$, we have the
following equality: 
\begin{equation*}
D_{\theta }(\overline{G}-G_{p}^{\Omega _{k}})=0.
\end{equation*}%
\noindent So $\overline{G}-G_{p}^{\Omega _{k}}$ is smooth near $p$ by the
regularity result for $\Delta _{b}$ and $D_{\theta }$ being "covariant" to $%
D_{\Theta }$ $=b_{n}\Delta _{b},$ $\Theta :$ standard contact form in the
Heisenberg group (cf. the argument in the end of the proof of Lemma 4.1).
Therefore we have the following result.

\bigskip

\textbf{Theorem 3.4. }\textit{Let }$G_{p}=\lim_{k\rightarrow \infty
}G_{p}^{\Omega _{k}}$\textit{. Then }$G_{p}$\textit{\ is a positive
fundamental solution of }$D_{\theta }$\textit{\ with pole at }$p.$\textit{\
Moreover, }$G_{p}$\textit{\ is minimal among all positive fundamental
solutions.}

\textit{\bigskip }

\proof
By the strong maximum principle of Bony and (\ref{barrier}), the sequence of
the Green's functions $\{G_{p}^{\Omega _{k}}\}$ is strictly increasing and
has an upper bound. Thus, away from $\Phi ^{-1}(y)$, the limit of $%
G_{p}^{\Omega _{k}}$ exists. Next by the standard argument using Bony's
maximum principle, $\Phi ^{-1}(y)\backslash \{p\}$ is a set of removable
singularities for $G_{p}$.

The minimality of $G_{p}$ follows from its construction, i.e., if $F_{p}$ is
another global positive fundamental solution on $\widetilde{M}$ with pole at 
$p$, then again Bony's maximum principle implies that $G_{p}^{\Omega
_{k}}<F_{p}$ for any $k$, and the conclusion follows.

\endproof%

\bigskip

We would like to discuss the decay properties of the constructed $G_{p}$, in
particular, its integrability away from the pole $p$. We define $s(M),$ the
minimum exponent of the integrability of $G_{p}$ by%
\begin{equation}
s(M):=\inf {\{s:\int_{\widetilde{M}\setminus U_{p}}G_{p}^{s}dV_{\theta
}<\infty \}}  \label{s(M)}
\end{equation}

\noindent where $U_{p}$ is a neighborhood of $p$.

\bigskip

\textbf{Theorem 3.5.} $s(M)$\textit{\ is a }$CR$\textit{\ invariant and
satisfies the following inequality: }%
\begin{equation}
s(M)\leq 1.  \label{3.12.1}
\end{equation}

\textit{\bigskip }

\proof
For $p$ $\in $ $\widetilde{M}$, let $U_{p}$ be a small neighborhood of $p$
with smooth boundary. Let $\{\Omega _{k}\subset \Omega _{k+1}:k=1,\cdots ,\}$
be a relatively compact, smooth exhaustion of the universal covering $%
\widetilde{M}$. We may assume that $p$ $\in $ $U_{p}$ $\subset $ $\Omega
_{1} $. Because the Dirichlet Green's function $G_{p}^{\Omega _{k}}$ is
smooth in $\Omega _{k-1}\backslash U_{p},$ Bony's maximum principle implies
that 
\begin{equation*}
\sup_{\Omega _{k-1}\backslash U_{p}}G_{p}^{\Omega _{k}}\leq \sup_{\Omega
_{k}\backslash U_{p}}G_{p}^{\Omega _{k}}\leq \sup_{\partial
U_{p}}G_{p}^{\Omega _{k}}<\max_{\partial U_{p}}G_{p}
\end{equation*}%
\noindent Therefore we may ($C^{\infty })$ smoothly extend $G_{p}^{\Omega
_{k}}$ into $U_{p}$ with the extension smaller than $\max_{\partial
U_{p}}G_{p},$ but positive. Denote this extension (which is a smooth
function over $\Omega _{k-1}$) by $\bar{G}_{p}^{\Omega _{k}}$.

For each $\alpha \geq 0$, let $u_{k}$ be the solution to the following
Dirichlet problem:%
\begin{equation*}
\begin{array}{rcl}
D_{\theta }(u) & = & (\bar{G}_{p}^{\Omega _{k}})^{\alpha }\ \text{in}\
\Omega _{k-1} \\ 
u|_{\partial \Omega _{k-1}} & = & 0.%
\end{array}%
\end{equation*}%
\noindent By the weak maximum principle, we obtain 
\begin{equation}
\max_{\Omega _{k-1}}u_{k}\leq \frac{(\max_{\partial U_{p}}G_{p})^{\alpha }}{%
R_{0}}  \label{3.13.1}
\end{equation}

\noindent where $R_{0}>0$ is the lower bound of the Tanaka-Webster scalar
curvature of $\widetilde{M}$. By the solution representation in $\Omega
_{k-1},$ we have

\begin{eqnarray}
u_{k}(p) &=&\int_{\Omega _{k-1}}G_{p}^{\Omega _{k-1}}(q)(\bar{G}_{p}^{\Omega
_{k}}(q))^{\alpha }dV_{\theta }(q)  \label{3.14} \\
&\geq &\int_{\Omega _{k-1}\backslash U_{p}}G_{p}^{\Omega _{k-1}}(q)(\bar{G}%
_{p}^{\Omega _{k}}(q))^{\alpha }dV_{\theta }(q)  \notag
\end{eqnarray}

\noindent By the monotonicity of $G_{p}^{\Omega _{k}}$ and letting $k$ $%
\rightarrow $ $\infty $, we conclude that 
\begin{equation*}
\int_{\widetilde{M}\setminus U_{p}}G_{p}^{1+\alpha }dV_{\theta }\leq \frac{%
(\max_{\partial U_{p}}G_{p})^{\alpha }}{R_{0}}.
\end{equation*}%
\noindent from (\ref{3.14}) and (\ref{3.13.1}). So (\ref{3.12.1}) follows.
Finally, we need to show that $s(M)$ is a well defined $CR$ invariant. It is
routine to check that the definition of $s(M)$ is independent of the choice
of $U_{p}$. Also $s(M)$ is independent of the choice of contact form from
the transformation law of Green's functions with respect to two different
contact forms (cf. Proposition 2.6 of Chapter VI in \cite{SY2} for the
conformal case).

\endproof%

\bigskip

Let $\rho $ denote the Carnot-Carath\'{e}odory distance on $\widetilde{M}$
with respect to the Levi metric (see, e.g., \cite{Str})$.$ Let $B_{r}(p)$ $%
\subset $ $\widetilde{M}$ denote the ball of radius $r$, centered at $p,$
with respect to the Carnot-Carath\'{e}odory distance $\rho .$ From Theorem
3.5 and a Moser's iteration procedure, we have the following result.

\bigskip

\textbf{Proposition 3.6.} \textit{There holds }%
\begin{equation*}
\lim_{r\rightarrow \infty }\left( \sup \{G_{p}(x):\rho (x,p)\geq r\}\right)
=0.
\end{equation*}

\textit{\bigskip }

\proof
Recall that $b_{n}$ $:=$ $2+\frac{2}{n}$. By Theorem 3.5,\ $\int_{\widetilde{%
M}\setminus U_{p}}G_{p}^{b_{n}}dV_{\theta }<\infty $. Thus we have 
\begin{equation*}
\lim_{r\rightarrow \infty }\int_{\{x:\rho (x,p)\geq
r\}}G_{p}^{b_{n}}dV_{\theta }=0.
\end{equation*}%
\noindent Therefore it is enough to establish the estimate 
\begin{equation*}
G_{p}(x)\leq C\left( \int_{B_{1}(x)}G_{p}^{b_{n}}dV_{\theta }\right)
^{1/b_{n}}\ \text{for all}\ x\in \widetilde{M}\backslash B_{2}(p)
\end{equation*}%
\noindent where $B_{1}(x)$ is a ball of radius $1,$ centered at $x$. First,
we have the equation for $G_{p}$ 
\begin{equation*}
\Delta _{b}G_{p}+\frac{1}{b_{n}}RG_{p}=0\ \ \text{on}\ \ \widetilde{M}%
\backslash B_{2}(p).
\end{equation*}%
Take $q\geq b_{n}$ $:=$ $2+\frac{2}{n}$. Multiplying the above formula by $%
G_{p}^{q-1}\phi ^{2}$ and integrating by parts give%
\begin{eqnarray}
&&(q-1)\int_{\widetilde{M}}\phi ^{2}G_{p}^{q-2}|\nabla
_{b}G_{p}|^{2}dV_{\theta }+\frac{1}{b_{n}}\int_{\widetilde{M}}RG_{p}^{q}\phi
^{2}dV_{\theta }  \label{3.15} \\
&\leq &2\int_{\widetilde{M}}\phi G_{p}^{q-1}|\nabla _{b}\phi ||\nabla
_{b}G_{p}|dV_{\theta }  \notag \\
&\leq &\alpha \int_{\widetilde{M}}\phi ^{2}G_{p}^{q-2}|\nabla
_{b}G_{p}|^{2}dV_{\theta }+\frac{1}{\alpha }\int_{\widetilde{M}}|\nabla
_{b}\phi |^{2}G_{p}^{q}dV_{\theta }  \notag
\end{eqnarray}

\noindent for all $\alpha >0$. Here $\phi \in C_{0}^{\infty }(\widetilde{M}%
\backslash B(p))$. Taking $\alpha =q-2$ in (\ref{3.15}), we get 
\begin{equation}
\int_{\widetilde{M}}\phi ^{2}G_{p}^{q-2}|\nabla _{b}G_{p}|^{2}dV_{\theta
}\leq \frac{1}{q-2}\int_{\widetilde{M}}|\nabla _{b}\phi
|^{2}G_{p}^{q}dV_{\theta }  \label{leq1}
\end{equation}%
\noindent by $R$ $>$ $0$ and $G_{p}$ $>$ $0.$ On the other hand, taking $%
\alpha $ $=$ $1$ in (\ref{3.15}) gives%
\begin{equation}
\begin{split}
2\int_{\widetilde{M}}\phi G_{p}^{q-1}|\nabla _{b}\phi ||\nabla
_{b}G_{p}|dV_{\theta }& \leq \int_{\widetilde{M}}\phi ^{2}G_{p}^{q-2}|\nabla
_{b}G_{p}|^{2}dV_{\theta }+\int_{\widetilde{M}}|\nabla _{b}\phi
|^{2}G_{p}^{q}dV_{\theta } \\
& \leq (\frac{q-1}{q-2})\int_{\widetilde{M}}|\nabla _{b}\phi
|^{2}G_{p}^{q}dV_{\theta }.
\end{split}
\label{leq2}
\end{equation}%
\noindent by (\ref{leq1}). It then follows from (\ref{leq2}) that 
\begin{equation}
\begin{split}
& \int_{\widetilde{M}}\left( |\nabla _{b}(\phi G_{p}^{q/2})|^{2}+\frac{1}{%
b_{n}}R\phi ^{2}G_{p}^{q}\right) dV_{\theta } \\
\leq & \int_{\widetilde{M}}\left( |\nabla _{b}\phi |^{2}G_{p}^{q}+\frac{q^{2}%
}{4}\phi ^{2}G_{p}^{q-2}|\nabla _{b}G_{p}|^{2}+q\phi G_{p}^{q-1}|\nabla
_{b}\phi ||\nabla _{b}G_{p}|+\frac{1}{b_{n}}R\phi ^{2}G_{p}^{q}\right)
dV_{\theta } \\
\leq & C_{n}q^{2}\int_{\widetilde{M}}\left( |\nabla _{b}\phi |^{2}+\phi
^{2}\right) G_{p}^{q}dV_{\theta }
\end{split}
\label{leq3}
\end{equation}%
\noindent for some constant $C_{n}$ independent of $q$. Applying the Sobolev
inequality 
\begin{equation*}
\left( \int_{\widetilde{M}}|\phi |^{b_{n}}dV_{\theta }\right) ^{2/b_{n}}\leq
\lambda (\widetilde{M})^{-1}\int_{\widetilde{M}}\left( |\nabla _{b}\phi
|^{2}+\frac{1}{b_{n}}R\phi ^{2}\right) dV_{\theta }
\end{equation*}%
\noindent (note that $\lambda (\widetilde{M})$ $>$ $0.$ In fact, $\lambda (%
\widetilde{M})$ $=$ $\lambda (S^{2n+1})$ a $CR$ analogue of Theorem 2.2 of
Chapter VI in \cite{SY2}), we obtain 
\begin{equation}
\begin{split}
& \left[ \int_{\widetilde{M}}(\phi G_{p}^{q/2})^{b_{n}}\right] ^{2/b_{n}} \\
\leq & \lambda (\widetilde{M})^{-1}\int_{\widetilde{M}}\left( |\nabla
_{b}(\phi G_{p}^{q/2})|^{2}+\frac{1}{b_{n}}R\phi ^{2}G_{p}^{q}\right)
dV_{\theta } \\
\leq & \tilde{C}_{n}q^{2}\int_{\widetilde{M}}\left( |\nabla _{b}\phi
|^{2}+\phi ^{2}\right) G_{p}^{q}dV_{\theta }
\end{split}
\label{leq4}
\end{equation}%
\noindent by (\ref{leq3}) for some constant $\tilde{C}_{n}$ independent of $%
q $.

We will use (\ref{leq4}) repeatedly with 
\begin{equation*}
q_{0}=b_{n}=2r,\ \ q_{k}=q_{0}r^{k},\ \ \text{with}\ \ r=\frac{n+1}{n}.
\end{equation*}%
\noindent Define a sequence of cut-off functions as follows. Set $%
a_{0}=1,a_{k}=1-\sum_{i=1}^{k}3^{-i}$ for $k\geq 1$, and we require that for
each $k$ the function $\phi _{k}\in C_{0}^{\infty }(\widetilde{M})$
satisfies 
\begin{equation*}
\phi _{k}=\{%
\begin{array}{rl}
1, & y\in B_{a_{k}}(x) \\ 
0, & y\not\in B_{a_{k-1}}(x),%
\end{array}%
\end{equation*}%
\begin{equation*}
0\leq \phi _{k}\leq 1,\ \ |\nabla _{b}\phi _{k}|\leq 2\cdot 3^{k}.
\end{equation*}%
\noindent Then iteratively we get from (\ref{leq4}) that

\begin{eqnarray}
&&\left( \int_{B_{a_{k}}(x)}G_{p}^{q_{k+1}}dV_{\theta }\right) ^{1/q_{k+1}}
\label{3.20} \\
&\leq &(Cq_{k}^{2})^{1/q_{k}}(4\cdot 3^{2k}+1)^{1/q_{k}}\left(
\int_{B_{a_{k-1}}(x)}G_{p}^{q_{k}}dV_{\theta }\right) ^{1/q_{k}}\leq ... 
\notag \\
&\leq &\prod_{j=1}^{k}(Cr^{2j})^{1/pr^{j}}\left(
\int_{B_{a_{0}}(x)}G_{p}^{b_{n}}dV_{\theta }\right) ^{1/b_{n}}  \notag
\end{eqnarray}
Since the product converges as $k\rightarrow \infty $, we can take the limit 
$k\rightarrow \infty $ in (\ref{3.20}) to get 
\begin{equation*}
\sup_{y\in B_{\frac{1}{2}}(x)}G_{p}(y)\leq C\left(
\int_{B_{1}(x)}G_{p}^{b_{n}}dV\right) ^{1/b_{n}}.
\end{equation*}%
This completes the proof.

\endproof%

\bigskip

\section{Injectivity of the $CR$ developing map}

Let $p$ $\in $ $\widetilde{M}$, and recall that $G_{p}$ denotes the minimal
positive Green's function for the $CR$ invariant sublaplacian $D_{\theta }$
with pole at $p$. Let $\Phi :\widetilde{M}\rightarrow S^{2n+1}$ be a $CR$
developing map. Recall that $H_{y}$ denotes the Green's function for the $CR$
invariant sublaplacian $D_{0}$ of $S^{2n+1}$ with the pole $y=\Phi (p)$. The
normalized pullback of $H_{y}$ is%
\begin{equation*}
\overline{G}:=|\Phi ^{^{\prime }}(p)|^{-(n+2)}|\Phi ^{^{\prime
}}|^{n}H_{y}\circ \Phi
\end{equation*}%
\noindent which has poles in $\Phi ^{-1}(y)$. Observe that $\Phi $ is one to
one if and only if $\Phi ^{-1}(y)=\{p\}$. Therefore to prove injectivity of $%
\Phi ,$ it suffices to prove $G_{p}=\overline{G}$.

For any $y\in S^{2n+1}$, the Cayley transform is a global $CR$
diffeomorphism 
\begin{equation*}
C_{y}:(S^{2n+1}\backslash \{y\},\theta _{S^{2n+1}})\rightarrow (H^{n},\Theta
)
\end{equation*}%
\noindent with $C_{y}(y)=\infty $ and 
\begin{equation*}
C_{y}^{\ast }\Theta =H_{y}^{\frac{2}{n}}\theta _{S^{2n+1}}.
\end{equation*}%
\noindent This means that $(S^{2n+1}\backslash \{y\},H_{y}^{\frac{2}{n}%
}\theta _{S^{2n+1}})$ is Heisenberg flat. Now we have

\begin{eqnarray}
|\Phi ^{^{\prime }}(p)|^{\frac{-2(n+2)}{n}}\Phi ^{\ast }\circ C_{y}^{\ast
}(\Theta ) &=&|\Phi ^{^{\prime }}(p)|^{\frac{-2(n+2)}{n}}\Phi ^{\ast
}(H_{y}^{\frac{2}{n}}\theta _{S^{2n+1}})  \label{4.1} \\
&=&|\Phi ^{^{\prime }}(p)|^{\frac{-2(n+2)}{n}}(H_{y}\circ \Phi )^{\frac{2}{n}%
}\Phi ^{\ast }\theta _{S^{2n+1}}  \notag \\
&=&|\Phi ^{^{\prime }}(p)|^{\frac{-2(n+2)}{n}}(H_{y}\circ \Phi )^{\frac{2}{n}%
}|\Phi ^{^{\prime }}|^{2}\theta  \notag \\
&=&(|\Phi ^{^{\prime }}(p)|^{-(n+2)}(H_{y}\circ \Phi )|\Phi ^{^{\prime
}}|^{n})^{\frac{2}{n}}\theta  \notag \\
&=&\overline{G}^{\frac{2}{n}}\theta .  \notag
\end{eqnarray}

\noindent It follows from, (\ref{4.1}) that $(\widetilde{M},$ $\bar{\theta}:=%
\overline{G}^{\frac{2}{n}}\theta )$ is Heisenberg flat away from $\Phi
^{-1}(y)$.

Cosider the quotient of $G_{p}$ and $\overline{G}:$%
\begin{equation*}
v:=\frac{G_{p}}{\overline{G}}.
\end{equation*}%
\noindent By (\ref{barrier}) we have $\overline{G}-G_{p}\geq 0$ away from $%
\Phi ^{-1}(y)$. So there holds 
\begin{equation*}
0<v\leq 1
\end{equation*}%
\noindent away from $\Phi ^{-1}(y).$ Taking $u$ $=$ $G_{p}$ and $\varphi $ $%
= $ $1$ in (\ref{TL}) , we obtain that on $\widetilde{M}\backslash \{p\}$: 
\begin{equation}
R(G_{p}^{\frac{2}{n}}\theta )=G_{p}^{-1-\frac{2}{n}}D_{\theta }(G_{p})=0.
\label{4.1.1}
\end{equation}%
\noindent Here for a contact form $\eta ,$ $R(\eta )$ or $R_{\eta }$ means
the Tanaka-Webster scalar curvature with respect to $\eta .$ Writing $G_{p}^{%
\frac{2}{n}}\theta =v^{\frac{2}{n}}\bar{\theta}$, we get 
\begin{equation*}
0=R(G_{p}^{\frac{2}{n}}\theta )=R(v^{\frac{2}{n}}\bar{\theta})=v^{-1-\frac{2%
}{n}}D_{\bar{\theta}}(v)
\end{equation*}%
\noindent away from $\Phi ^{-1}(y)$ by (\ref{4.1.1}) and (\ref{TL}) with $u$ 
$=$ $v$ and $\varphi $ $=$ $1.$ Therefore we have 
\begin{equation}
(b_{n}\bar{\Delta}_{b}+R_{\bar{\theta}})(v)=b_{n}\bar{\Delta}_{b}(v)=0
\label{harofv}
\end{equation}

\noindent by noting that $R_{\bar{\theta}}=0.$ We would like to examine the
asymptotic behavior of $v$ near $\Phi ^{-1}(y)$. We will often write the
coordinates $(z_{1},$ $...,$ $z_{n},$ $t)$ in $H^{n}$ by $(z,$ $t)$ where $z$
$=$ $(z_{1},$ $...,$ $z_{n}).$ Define the Heisenberg norm $\rho _{0}(z,t)$ $%
:=$ $(|z|^{4}+t^{2})^{1/4}.$ Denote the fundamental solution to $D_{\Theta }$
$=$ $b_{n}\Delta _{b}$ by $c(n)\rho _{0}(z,t)^{-2n}$ for some dimensional
constant $c(n)$ (\cite{FS}).

\bigskip

\textbf{Lemma 4.1.} \textit{For each }$\bar{p}\in \Phi ^{-1}(y),$\textit{\
we can choose a coordinate map }$(z,t):\widetilde{M}\rightarrow H^{n}$%
\textit{\ and a smooth function }$u$ \textit{\ near }$\bar{p}$\textit{\ such
that }$(z(\bar{p}),t(\bar{p}))$\textit{\ }$=$\textit{\ }$(0,0)$\textit{\ and
there hold}%
\begin{equation}
G_{p}(z,t)=c(n)u(p)u\rho _{0}(z,t)^{-2n}+a\text{ }smooth\text{ }function
\label{4.3.1}
\end{equation}%
\textit{\noindent near }$\bar{p}$ $=$ $p$ \textit{\ and}%
\begin{equation}
\overline{G}(z,t)=c(n)u(q)u\rho _{0}(z,t)^{-2n}+a\text{ }smooth\text{ }%
function  \label{4.3.2}
\end{equation}%
\textit{\noindent near }$\bar{p}$ $=$ $q$,\textit{\ resp..}

\bigskip

\proof
Let $-y\in S^{2n+1}$ be the antipodal point of $y\in S^{2n+1}$. Consider the
Cayley transform $C_{-y}$ (with pole at $-y$). Obviously, $C_{-y}(y)=0$ $\in 
$ $H^{n}$ and 
\begin{equation}
(C_{-y}^{-1})^{\ast }(H_{-y}^{\frac{2}{n}}\theta _{S^{2n+1}})=\Theta
\label{4.3.3}
\end{equation}%
\noindent It follows from (\ref{4.3.3}) that 
\begin{equation}
(C_{-y}\circ \Phi )^{\ast }\Theta =(H_{-y}\circ \Phi )^{\frac{2}{n}}|\Phi
^{^{\prime }}|^{2}\theta .  \label{4.3.5}
\end{equation}%
\noindent Here we have written $\Phi ^{\ast }(\theta _{S^{2n+1}})=|\Phi
^{^{\prime }}|^{2}\theta .$ Let $u$ :$=$ $(H_{-y}\circ \Phi )|\Phi
^{^{\prime }}|^{n}.$ We can then write (\ref{4.3.5}) as%
\begin{equation}
(C_{-y}\circ \Phi )^{\ast }\Theta =u^{\frac{2}{n}}\theta  \label{4.3.6}
\end{equation}

\noindent near $q$ $\in $ $\Phi ^{-1}(y).$ Take $C_{-y}\circ \Phi $ $:$ $%
\widetilde{M}\rightarrow H^{n}$ as a coordinate map $(z,t)$. By (\ref{4.3.6}%
) and (\ref{TL}) in Section 2 with $\varphi $ $=$ $c(n)\rho _{0}(z,t)^{-2n}$%
, we obtain%
\begin{eqnarray*}
D_{\theta }(uc(n)\rho _{0}(z,t)^{-2n}) &=&u^{1+\frac{2}{n}}D_{\Theta
}(c(n)\rho _{0}(z,t)^{-2n}) \\
&=&u^{1+\frac{2}{n}}\delta _{(0,0)}.
\end{eqnarray*}%
\noindent Note that the volume change formula is ($(C_{-y}\circ \Phi )^{\ast
})$ $\Theta \wedge (d\Theta )^{n}$ $=$ $u^{2+\frac{2}{n}}\theta \wedge
(d\theta )^{n}.$ So in view of Theorem 3.4, (\ref{3.11.1}), and (\ref{3.11.2}%
),.we have%
\begin{eqnarray}
D_{\theta }(G_{p}-u(p)uc(n)\rho _{0}(z,t)^{-2n}) &=&0,  \label{4.3.7} \\
D_{\theta }(\bar{G}-u(\bar{p})uc(n)\rho _{0}(z,t)^{-2n}) &=&0  \notag
\end{eqnarray}%
\noindent near $p,$ $\bar{p}$ $\in $ $\Phi ^{-1}(y)$, resp.. By applying (%
\ref{TL}) we obtain%
\begin{eqnarray}
D_{(C_{-y}\circ \Phi )^{\ast }\Theta }(u^{-1}(G_{p}-u(p)uc(n)\rho
_{0}(z,t)^{-2n}) &=&0  \label{4.3.8} \\
D_{(C_{-y}\circ \Phi )^{\ast }\Theta }(u^{-1}(\bar{G}-u(\bar{p})uc(n)\rho
_{0}(z,t)^{-2n}) &=&0  \notag
\end{eqnarray}%
\noindent according to (\ref{4.3.6}) and (\ref{4.3.7}). Observe that $%
R_{(C_{-y}\circ \Phi )^{\ast }\Theta }$ $=$ $(C_{-y}\circ \Phi )^{\ast
}R_{\Theta }$ $=$ $R_{\Theta }\circ (C_{-y}\circ \Phi )$ $=$ $0.$ Hence $%
D_{(C_{-y}\circ \Phi )^{\ast }\Theta }$ $=$ $b_{n}\Delta _{b}.$ By the
regularity result for $\Delta _{b}$ (e.g., Theorem 16.7 in \cite{FS}), we
have (\ref{4.3.1}) and (\ref{4.3.2}) from (\ref{4.3.8}) and $G_{p},$.$\bar{G}%
,$ $\rho _{0}(z,t)^{-2n}$ being $L_{loc}^{1}.$

\endproof%

\bigskip

By Lemma 4.1 we deduce that near $p$ (which is a pole of $\overline{G}$): 
\begin{equation*}
v(z,t)=1+O(\rho _{0}^{2n}),
\end{equation*}%
\noindent and near $\bar{p}$ $\in $ $\Phi ^{-1}(y)\setminus \{p\}$ 
\begin{equation*}
v(z,t)=O(\rho _{0}^{2n}).
\end{equation*}%
\noindent since $G_{p}$ is smooth near any $\bar{p}$ $\in $ $\Phi
^{-1}(y)\setminus \{p\}$. From Lemma 4.1 we also have

\begin{eqnarray*}
|\nabla _{b}v| &=&O(\rho _{0}^{2n-1}), \\
|\nabla _{b}|\nabla _{b}v|| &=&O(\rho _{0}^{2n-2}),\Delta _{b}v=O(\rho
_{0}^{2n-2})
\end{eqnarray*}

\noindent near any $\bar{p}$ $\in $ $\Phi ^{-1}(y).$

With the asymptotic behavior of $G_{p},$ $\bar{G},$ and $v$ near any $\bar{p}
$ $\in $ $\Phi ^{-1}(y)$ understood, we observe that the set $\Phi ^{-1}(y)$
has no contribution when we play the integration by parts in the computation
below throughout the remaining section. We would like to show that $v$ is a
constant, and hence is identically one since $v(p)=1$. Write $G,$ $dV$
instead of $G_{p},$ $dV_{\theta }$ for short in the remaining section. Also
note that the notation $C$ may mean different constants.

\bigskip

\textbf{Lemma 4.2.} \textit{There exists a constant }$C>0$\textit{\ such
that for any }$\phi \in C_{0}^{\infty }(\widetilde{M})$\textit{\ there holds}%
\begin{equation}
\int_{\widetilde{M}}\phi ^{2}|\nabla _{b}\log {\overline{G}}|^{2}dV\leq
C\int_{\widetilde{M}}(\phi ^{2}|\nabla _{b}\log {G}|^{2}+|\nabla _{b}\phi
|^{2})dV.  \label{firequ}
\end{equation}

\bigskip

\proof
Away from $\Phi ^{-1}(y)$ we have 
\begin{equation}
\begin{split}
\Delta _{b}\log {\overline{G}}& =\overline{G}^{-1}\Delta _{b}\overline{G}-(-%
\overline{G}^{-2}|\nabla _{b}\overline{G}|^{2}) \\
& =\overline{G}^{-1}\Delta _{b}\overline{G}+|\nabla _{b}\log {\overline{G}}%
|^{2} \\
& =-b_{n}^{-1}R(\theta )+|\nabla _{b}\log {\overline{G}}|^{2} \\
& =G^{-1}\Delta _{b}G+|\nabla _{b}\log {\overline{G}}|^{2} \\
& =\Delta _{b}\log {G}-|\nabla _{b}\log {G}|^{2}+|\nabla _{b}\log {\overline{%
G}}|^{2}
\end{split}
\label{4.3.4}
\end{equation}

\noindent Multiplying $\phi ^{2}$ on both sides of (\ref{4.3.4}) and
integrating by parts, we obtain 
\begin{equation*}
\begin{split}
\int_{\widetilde{M}}\phi ^{2}|\nabla _{b}\log {\overline{G}}|^{2}dV& =\int_{%
\widetilde{M}}\phi ^{2}|\nabla _{b}\log {G}|^{2}dV+2\int_{\widetilde{M}}\phi
\nabla _{b}\phi (\nabla _{b}\log {\overline{G}}-\nabla _{b}\log {G})dV \\
& \leq C\int_{\widetilde{M}}(\phi ^{2}|\nabla _{b}\log {G}|^{2}+|\nabla
_{b}\phi |^{2})dV.
\end{split}%
\end{equation*}%
\noindent for some constant $C$ by noting that the boundary integral around $%
\bar{p}$ $\in $ $\Phi ^{-1}(y)$ tends to zero$.$ Here we have used the
Schwarz inequality with $\varepsilon $ in the last inequality.

\endproof%

\bigskip

Let $\bar{T}$ $(\overline{\nabla }_{b}, \overline{\Delta }_{b}$, $d\overline{%
V},$ $\bar{R}$ $=$ $R_{\bar{\theta}},$ etc., resp.) denote the corresponding
quantity of $T$ ($\nabla _{b},$ $\Delta _{b},$ $dV,$ $R,$ etc., resp.) with
respect to $\bar{\theta}$ (while fixing $J$)$.$

\bigskip

\textbf{Lemma 4.3.} \textit{There holds }$v_{\bar{0}}$\textit{\ }$:=$\textit{%
\ }$\overline{T}v\equiv 0$\textit{\ in either of the following cases:}

\textit{(a) }$n\geq 3$

\textit{(b) }$n=2$\textit{\ and }$s(M)<1.$

\bigskip

\proof
First observe that the Paneitz-like operator $P$ is nonnegative for $\varphi 
$ $\in $ $C_{0}^{\infty }(\widetilde{M})$ in \ref{NG} if $n$ $\geq $ $2$
(Extending Theorem 3.2 in \cite{CC} to this situation). With respect to $%
\bar{\theta}$ (Heisenberg flat)$,$ the torsion vanishes and hence $\kappa $ $%
=$ $0$ in (\ref{IBF2}). Therefore by extending (\ref{IBF2}) to the situation
for $\bar{\theta}$ (singular at $\bar{p}$ $\in $ $\Phi ^{-1}(y))$ in view of
the asymptotic behavior of $v$ discussed before$,$ we have 
\begin{equation}
\begin{split}
n^{2}\int_{\widetilde{M}}(\phi v)_{\bar{0}}^{2}d\bar{V}& \leq \int_{%
\widetilde{M}}(\overline{\Delta }_{b}(\phi v))^{2}d\overline{V} \\
& =\int_{\widetilde{M}}(\phi \overline{\Delta }_{b}v+v\overline{\Delta }%
_{b}\phi -2\overline{\nabla }_{b}\phi \overline{\nabla }_{b}v)^{2}d\overline{%
V} \\
& \leq C\left( \int_{\widetilde{M}}v^{2}(\overline{\Delta }_{b}\phi )^{2}d%
\overline{V}+\int_{\widetilde{M}}|\overline{\nabla }_{b}\phi |^{2}|\overline{%
\nabla }_{b}v|^{2}d\overline{V}\right) \\
& =C(I+II)
\end{split}
\label{24}
\end{equation}%
\noindent for $\phi $ $\in $ $C_{0}^{\infty }(\widetilde{M}),$ where $%
I=\int_{\widetilde{M}}v^{2}(\overline{\Delta }_{b}\phi )^{2}d\overline{V}$
and $II=\int_{\widetilde{M}}|\overline{\nabla }_{b}\phi |^{2}|\overline{%
\nabla }_{b}v|^{2}d\overline{V}$. Rewrite%
\begin{eqnarray}
II &=&\int_{\widetilde{M}}|\overline{\nabla }_{b}\phi |^{2}|\overline{\nabla 
}_{b}v|^{2}d\overline{V}  \label{24.0} \\
&=&\int_{\widetilde{M}}|\nabla _{b}\phi |^{2}|\nabla _{b}v|^{2}\overline{G}^{%
\frac{2n-2}{n}}dV  \notag
\end{eqnarray}

\noindent where%
\begin{eqnarray}
|\nabla _{b}v|^{2} &=&|\frac{\overline{G}\nabla _{b}G-G\nabla _{b}\overline{G%
}}{\overline{G}^{2}}|^{2}  \label{24.1} \\
&=&|v\nabla _{b}\ log\ G-v\nabla _{b}\log \overline{G}|^{2}  \notag \\
&\leq &C(v^{2}|\nabla _{b}\log {G}|^{2}+v^{2}|\nabla _{b}\log {\overline{G}}%
|^{2}).  \notag
\end{eqnarray}%
\noindent Let $q^{\prime }$ :$=$ $\frac{2(n-1)}{n}$. Then $q^{\prime }>1$ if
and only if $n\geq 3$. From (\ref{24.1}) we have%
\begin{equation}
\begin{split}
\overline{G}^{q^{\prime }}|\nabla _{b}v|^{2}& \leq C(v^{2}\overline{G}%
^{q^{\prime }}|\nabla _{b}\log {G}|^{2}+v^{2}\overline{G}^{q^{\prime
}}|\nabla _{b}\log {\overline{G}}|^{2}) \\
& =Cv^{2-q^{\prime }}(G^{q^{\prime }}|\nabla _{b}\log {G}|^{2}+G^{q^{\prime
}}|\nabla _{b}\log {\overline{G}}|^{2}) \\
& \leq CG^{q^{\prime }}(|\nabla _{b}\log {G}|^{2}+|\nabla _{b}\log {%
\overline{G}}|^{2}).
\end{split}
\label{25}
\end{equation}%
\noindent Recall that $B_{\rho }(p)$ $\subset $ $\widetilde{M}$ denote the
ball of radius $\rho $, centered at $p,$ with respect to the Carnot-Carath%
\'{e}odory distance. Substituting (\ref{25}) into (\ref{24.0}), we get%
\begin{equation}
II\leq C\rho ^{-2}\int_{B_{2\rho }(p)\backslash B_{\rho }(p)}G^{q^{\prime
}}(|\nabla _{b}\log {G}|^{2}+|\nabla _{b}\log {\overline{G}}|^{2})dV
\label{25.1}
\end{equation}

\noindent by taking a cutoff function $\phi $ such that $0$ $\leq $ $\phi $ $%
\leq $ $1,$ $\phi $ $=$ $1$ on $B_{\rho }(p),$ $\phi $ $=$ $0$ on ($%
\widetilde{M}\TEXTsymbol{\backslash}$$B_{2\rho }(p)),$ and $|\nabla _{b}\phi
| $ $\leq $ $\frac{C}{\rho }.$

Taking $\phi =\psi G^{\frac{q^{\prime }}{2}}$ in (\ref{firequ}), $\psi \in
C_{0}^{\infty }(\widetilde{M}\backslash \{p\}),$ we get 
\begin{equation}
\begin{split}
& \int_{\widetilde{M}}\psi ^{2}G^{q^{\prime }}|\nabla _{b}\log {\overline{G}}%
|^{2}dV \\
& \leq C\int_{\widetilde{M}}(\psi ^{2}G^{q^{\prime }}|\nabla _{b}\log {G}%
|^{2}+|\nabla _{b}(G^{q^{\prime }/2}\psi )|^{2})dV \\
& \leq C\int_{\widetilde{M}}(\psi ^{2}|\nabla _{b}\log {G}|^{2}+|\nabla
_{b}\psi |^{2})G^{q^{\prime }}dV.
\end{split}
\label{26}
\end{equation}%
\noindent Note that the integral of $|\nabla _{b}\log {G}|^{2}G^{q^{\prime
}} $ over a region containing $p$ diverges (this is why we need $\psi $
compactly supported away from $p).$ Choosing $\psi $ such that $0$ $\leq $ $%
\psi $ $\leq $ $1,\psi $ $=$ $1$ on $B_{2\rho }(p)\backslash B_{\rho }(p),$ $%
\psi $ $=$ $0$ on $B_{\rho /2}(p)$ $\cup $ ($\widetilde{M}\TEXTsymbol{%
\backslash}$$B_{4\rho }(p))$, and $|\nabla _{b}\psi |$ $\leq $ $\frac{C}{%
\rho },$ we get%
\begin{equation}
II\leq C\rho ^{-2}\int_{B_{4\rho }(p)\backslash B_{\rho /2}(p)}G^{q^{\prime
}}(1+|\nabla _{b}\log {G}|^{2})dV  \label{26.1}
\end{equation}

\noindent from (\ref{25.1}) and (\ref{26}) (for $\rho $ large).

For $I=\int_{\widetilde{M}}v^{2}(\overline{\Delta }_{b}\phi )^{2}d\overline{V%
}$, since 
\begin{equation*}
(b_{n}\overline{\Delta }_{b}+\overline{R})\phi =\overline{G}^{-1-\frac{2}{n}%
}(b_{n}\Delta _{b}+R_{\theta })(\overline{G}\phi ),
\end{equation*}%
\noindent we have 
\begin{equation}
\begin{split}
b_{n}\overline{\Delta }_{b}\phi & =\overline{G}^{-1-\frac{2}{n}}(b_{n}(\phi
\Delta _{b}\overline{G}+\overline{G}\Delta _{b}\phi -2\nabla _{b}\phi \nabla
_{b}\overline{G})+R_{\theta }\overline{G}\phi ) \\
& =\overline{G}^{-1-\frac{2}{n}}(\phi (b_{n}\Delta _{b}\overline{G}%
+R_{\theta }\overline{G})+b_{n}(\overline{G}\Delta _{b}\phi -2\nabla
_{b}\phi \nabla _{b}\overline{G})) \\
& =b_{n}(\overline{G}^{\frac{-2}{n}}\Delta _{b}\phi -2\overline{G}^{-1-\frac{%
2}{n}}\nabla _{b}\phi \nabla _{b}\overline{G}),
\end{split}
\label{27}
\end{equation}%
\noindent that is, 
\begin{equation}
\overline{\Delta }_{b}\phi =\overline{G}^{\frac{-2}{n}}\Delta _{b}\phi -2%
\overline{G}^{-1-\frac{2}{n}}\nabla _{b}\phi \nabla _{b}\overline{G}.
\label{27.1}
\end{equation}%
\noindent By (\ref{27.1}) we have%
\begin{eqnarray}
I &=&\int_{\widetilde{M}}v^{2}(\overline{G}^{\frac{-2}{n}}\Delta _{b}\phi -2%
\overline{G}^{-1-\frac{2}{n}}\nabla _{b}\phi \nabla _{b}\overline{G})^{2}d%
\overline{V}  \label{27.2} \\
&\leq &C(\int_{\widetilde{M}}v^{2}(\frac{\Delta _{b}\phi }{\overline{G}^{%
\frac{2}{n}}})^{2}d\overline{V}+\int_{\widetilde{M}}v^{2}\overline{G}^{-2-%
\frac{4}{n}}|\nabla _{b}\phi |^{2}|\nabla _{b}\overline{G}|^{2}d\overline{V})
\notag \\
&=&C(III+IV)  \notag
\end{eqnarray}

\noindent where%
\begin{eqnarray}
III &=&\int_{\widetilde{M}}v^{2}(\frac{\Delta _{b}\phi }{\overline{G}^{\frac{%
2}{n}}})^{2}d\overline{V}  \label{27.3} \\
&=&\int_{\widetilde{M}}(\frac{\Delta _{b}\phi }{\overline{G}^{\frac{2}{n}}}%
)^{2}(\frac{G}{\overline{G}})^{2}(\overline{G}^{\frac{2}{n}})^{n+1}dV  \notag
\\
&=&\int_{\widetilde{M}}(\Delta _{b}\phi )^{2}G^{2}\overline{G}^{\frac{2(n+1)%
}{n}-\frac{4}{n}-2}dV  \notag \\
&\leq &C\rho ^{-2}\int_{\widetilde{M}\backslash B_{\rho /4}}G^{2-\frac{2}{n}%
}dV  \notag
\end{eqnarray}

\noindent and%
\begin{eqnarray}
IV &=&\int_{\widetilde{M}}v^{2}\overline{G}^{-2-\frac{4}{n}}|\nabla _{b}\phi
|^{2}|\nabla _{b}\overline{G}|^{2}dV  \label{27.4} \\
&=&\int_{\widetilde{M}}v^{2}\overline{G}^{\frac{-2}{n}}|\nabla _{b}\phi
|^{2}|\nabla _{b}\overline{G}|^{2}dV  \notag \\
&=&\int_{\widetilde{M}}v^{2}\overline{G}^{2-\frac{2}{n}}|\nabla _{b}\phi
|^{2}|\nabla _{b}\log {\overline{G}}|^{2}dV  \notag \\
&\leq &\int_{\widetilde{M}}|\nabla _{b}\phi |^{2}G^{q^{\prime }}|\nabla
_{b}\log {\overline{G}}|^{2}dV  \notag \\
&\leq &C\rho ^{-2}\int_{B_{2\rho }(p)\backslash B_{\rho }(p)}G^{q^{\prime
}}|\nabla _{b}\log {\overline{G}}|^{2}dV  \notag
\end{eqnarray}

\noindent by a suitable choice of $\phi .$ Again we let $\phi =\psi
G^{q^{\prime }/2}$ in Lemma 4.2, where $\psi \in C_{0}^{\infty }(\widetilde{M%
}\setminus \{p\})$. By choosing $\psi $ suitably, we can convert (\ref{27.4}%
) into%
\begin{equation}
IV\leq C\rho ^{-2}\int_{B_{4\rho }(p)\backslash B_{\rho /2}(p)}G^{q^{\prime
}}(1+|\nabla _{b}\log {G}|^{2})dV.  \label{27.5}
\end{equation}%
Finally we are going to show that 
\begin{equation}
II\text{ (}IV,\text{ resp.)}\leq C\rho ^{-2}\int_{\widetilde{M}\backslash
B_{\rho /4}}G^{q^{\prime }}dV  \label{27.6}
\end{equation}%
\noindent (recall that $q^{\prime }$ $=$ $2$ $-$ $\frac{2}{n})$ for $n$ $%
\geq $ $3,$ and 
\begin{equation}
II\text{ (}IV,\text{ resp.)}\leq C\rho ^{-2}\int_{\widetilde{M}\backslash
B_{\rho /4}}G^{\widetilde{q}}dV  \label{27.6b}
\end{equation}

\noindent for $n$ $=$ $2,$ in which $\widetilde{q}$ $<$ $1.$ First consider
the case $(a)$ $n\geq 3$. So $q^{\prime }>1$. Now $b_{n}\Delta
_{b}G+R_{\theta }G$ $=$ $0$ (away from $p)$ and $R_{\theta }$ $>$ $C$ $>$ $0$%
. This implies that 
\begin{equation*}
\Delta _{b}G=-b_{n}^{-1}R_{\theta }G\leq -b_{n}^{-1}CG
\end{equation*}%
\noindent Multiplying by $\phi ^{2}G^{q^{\prime }-1}$ with $\phi \in
C_{0}^{\infty }(\widetilde{M}\setminus \{p\})$ and integrating by parts give

\begin{eqnarray*}
0 &\geq &\int_{\widetilde{M}}\phi ^{2}G^{q^{\prime }-1}\Delta _{b}G\ dV \\
&=&\int_{\widetilde{M}}\nabla _{b}(\phi ^{2}G^{q^{\prime }-1})\nabla _{b}GdV
\\
&=&\int_{\widetilde{M}}2\phi G^{q^{\prime }-1}\nabla _{b}\phi \nabla
_{b}GdV+(q^{\prime }-1)\int_{\widetilde{M}}\phi ^{2}G^{q^{\prime }-2}|\nabla
_{b}G|^{2}dV.
\end{eqnarray*}%
\noindent By Young's inequality with $\varepsilon $, we get 
\begin{equation}
\int_{\widetilde{M}}\phi ^{2}G^{q^{\prime }-2}|\nabla _{b}G|^{2}dV\leq
C\int_{\widetilde{M}}G^{q^{\prime }}|\nabla _{b}\phi |^{2}dV  \label{27.7}
\end{equation}

\noindent (noting that we have used the fact $q^{\prime }>1).$ By choosing $%
\phi $ appropriately in (\ref{27.7}) and observing that $G^{q^{\prime
}-2}|\nabla _{b}G|^{2}$ $=$ $G^{q^{\prime }}|\nabla _{b}\log G|^{2}$, we
obtain%
\begin{equation}
\int_{B_{4\rho }(p)\backslash B_{\rho /2}(p)}G^{q^{\prime }}|\nabla _{b}\log
G|^{2}\leq C\rho ^{-2}\int_{B_{8\rho }(p)\backslash B_{\rho
/4}(p)}G^{q^{\prime }}dV.  \label{27.8}
\end{equation}

\noindent So we have (\ref{27.6}) for $II$ ($IV,$ resp.) by (\ref{26.1}) and
(\ref{27.8}) ((\ref{27.5}) and (\ref{27.8}), resp.).

For the case $(b)$, $n$ $=$ $2$ implies $q^{\prime }$ $=$ $1$. From $%
b_{n}\Delta _{b}G+R_{\theta }G=0$ where $|R_{\theta }|\leq c$ we have

\begin{equation*}
\Delta _{b}G=-b_{n}^{-1}R_{\theta }G\geq -b_{n}^{-1}cG,
\end{equation*}%
\noindent that is, 
\begin{equation*}
0\leq \Delta _{b}G+b_{n}^{-1}cG
\end{equation*}

\noindent Multiplying by $\phi ^{2}G^{\widetilde{q}-1}$ with $\widetilde{q}%
<1,\ \phi $ $\in $ $C_{0}^{\infty }(M\backslash \{p\})$ and integrating by
parts give%
\begin{equation}
\begin{split}
0& \leq \int_{\widetilde{M}}\phi ^{2}G^{\widetilde{q}-1}\Delta
_{b}GdV+b_{n}^{-1}c\int_{\widetilde{M}}\phi ^{2}G^{\widetilde{q}}dV \\
& =\int_{\widetilde{M}}\nabla _{b}(\phi ^{2}G^{\widetilde{q}-1})\nabla
_{b}GdV+b_{n}^{-1}c\int_{\widetilde{M}}\phi ^{2}G^{\widetilde{q}}dV \\
& =\int_{\widetilde{M}}2\phi G^{\widetilde{q}-1}\nabla _{b}\phi \nabla
_{b}GdV+b_{n}^{-1}c\int_{\widetilde{M}}\phi ^{2}G^{\widetilde{q}}dV+(%
\widetilde{q}-1)\int_{\widetilde{M}}\phi ^{2}G^{\widetilde{q}-2}|\nabla
_{b}G|^{2}dV.
\end{split}
\label{29}
\end{equation}%
\noindent From (\ref{29}) we have%
\begin{eqnarray*}
&&(1-\widetilde{q})\int_{\widetilde{M}}\phi ^{2}G^{\widetilde{q}-2}|\nabla
_{b}G|^{2}dV \\
&\leq &2\int_{\widetilde{M}}\phi G^{\widetilde{q}-1}|\nabla _{b}\phi
||\nabla _{b}G|dV+b_{n}^{-1}c\int_{\widetilde{M}}\phi ^{2}G^{\widetilde{q}}dV
\end{eqnarray*}

\noindent By Young's inequality with $\varepsilon $, we obtain%
\begin{equation*}
\int_{\widetilde{M}}\phi ^{2}G^{\widetilde{q}-2}|\nabla _{b}G|^{2}dV\leq
C(\int_{\widetilde{M}}G^{\widetilde{q}}|\nabla _{b}\phi |^{2}dV+\int_{%
\widetilde{M}}G^{\widetilde{q}}\phi ^{2}dV).
\end{equation*}%
\noindent Since $G^{\widetilde{q}-2}\geq G^{-1}$ on $\widetilde{M}\setminus
K $ for some compact subset $K$ by Proposition 3.6, we have 
\begin{equation}
\int_{\widetilde{M}}\phi ^{2}G^{-1}|\nabla _{b}G|^{2}dV\leq C(\int_{%
\widetilde{M}}G^{\widetilde{q}}|\nabla _{b}\phi |^{2}dV+\int_{\widetilde{M}%
}G^{\widetilde{q}}\phi ^{2}dV).  \label{29.1}
\end{equation}

\noindent Observing that $G^{-1}|\nabla _{b}G|^{2}$ $=$ $G|\nabla _{b}\log
G|^{2}$ and choosing a suitable cutoff function $\phi $ in (\ref{29.1})$,$
we obtain%
\begin{equation}
\int_{B_{4\rho }(p)\backslash B_{\rho /2}(p)}G|\nabla _{b}\log G|^{2}dV\leq
C(\frac{1}{\rho ^{2}}+1)\int_{B_{8\rho }(p)\backslash B_{\rho /4}(p)}G^{%
\widetilde{q}}dV.  \label{29.2}
\end{equation}

\noindent Thus we have (\ref{27.6b}) for $II$ ($IV,$ resp.) by (\ref{26.1})
and (\ref{29.2}) ((\ref{27.5}) and (\ref{29.2}), resp.) for $\rho $ large in
view of Proposition 3.6.

By (\ref{3.12.1}) and the assumption $s(M)$ $<$ $1$ for $n$ $=$ $2,$ we have
the convergence of the integrals in (\ref{27.6}) and (\ref{27.6b}) (in fact,
both of them tend to zero as $\rho $ $\rightarrow $ $\infty ).$ So as $\rho $
$\rightarrow $ $\infty ,$ $II$ and $IV$ go to zero. On the other hand, it is
clear that $III$ goes to zero as $\rho $ $\rightarrow $ $\infty $ by (\ref%
{27.3}) and (\ref{3.12.1}). So from (\ref{24}) and (\ref{27.2}) we conclude
that $v_{\bar{0}}$ $=$ $0.$

\endproof%

\bigskip

\proof
\textbf{(of Theorem A)}

We need only to prove $v$ $:=\frac{G}{\overline{G}}\equiv 1$. Let $q$ $=$ $%
\frac{2n}{n+1}$. We first prove that for any $\phi \in C_{0}^{\infty }(%
\widetilde{M})$ there holds%
\begin{equation}
\int_{\widetilde{M}}\phi ^{2}|\overline{\nabla }_{b}v|^{q-2}|\overline{%
\nabla }_{b}|\overline{\nabla }_{b}v||^{2}d\overline{V}\leq C\int_{%
\widetilde{M}}|\nabla _{b}\phi |^{2}\overline{G}^{q}|\nabla _{b}v|^{q}dV
\label{8}
\end{equation}%
\noindent for some constant $C$. Since $\overline{\Delta }_{b}v=0,$ $\bar{%
\theta}=\overline{G}^{\frac{2}{n}}\theta $ is flat (hence $Ric$ and $Tor$
vanish), and $v_{\bar{0}}$ $=$ $0$ by Lemma 4.3, we reduce the Bochner
formula (\ref{BF}) to

\begin{equation}
\frac{1}{2}\overline{\Delta }_{b}|\overline{\nabla }_{b}v|^{2}=-|(\overline{%
\nabla }^{H})^{2}v|^{2}.  \label{2}
\end{equation}%
\noindent Observe that 
\begin{equation}
|\overline{\nabla }_{b}|\overline{\nabla }_{b}v||^{2}\leq |(\overline{\nabla 
}^{H})^{2}v|^{2}.  \label{3}
\end{equation}%
\noindent For $q>1$ we compute 
\begin{equation}
\begin{split}
\overline{\Delta }_{b}|\overline{\nabla }_{b}v|^{q}& =\overline{\Delta }%
_{b}(|\overline{\nabla }_{b}v|^{2})^{\frac{q}{2}} \\
& =\frac{q}{2}(|\overline{\nabla }_{b}v|^{2})^{\frac{q}{2}-1}\overline{%
\Delta }_{b}|\overline{\nabla }_{b}v|^{2}-\frac{q}{2}(\frac{q}{2}-1)(|%
\overline{\nabla }_{b}v|^{2})^{\frac{q}{2}-2}|\overline{\nabla }_{b}|%
\overline{\nabla }_{b}v|^{2}|^{2} \\
& =\frac{q}{2}|\overline{\nabla }_{b}v|^{q-2}\overline{\Delta }_{b}|%
\overline{\nabla }_{b}v|^{2}-\frac{q(q-2)}{4}|\overline{\nabla }_{b}v|^{q-4}|%
\overline{\nabla }_{b}|\overline{\nabla }_{b}v|^{2}|^{2} \\
& \leq -q(q-1)|\overline{\nabla }_{b}v|^{q-2}|\overline{\nabla }_{b}|%
\overline{\nabla }_{b}v||^{2} \\
& =-C_{q}|\overline{\nabla }_{b}v|^{q-2}|\overline{\nabla }_{b}|\overline{%
\nabla }_{b}v||^{2}
\end{split}
\label{4}
\end{equation}%
\noindent where $C_{q}$ $=$ $q(q-1)$ $>$ $0$ for $n$ $\geq $ $2.$ For the
inequality in (\ref{4}) we have used (\ref{2}), (\ref{3}). Consider first
the case $\phi \in C_{0}^{\infty }(\widetilde{M}\setminus \Phi ^{-1}(y))$.
Multiplying (\ref{4}) by $\phi ^{2}$ and integrating by parts, we get

\begin{equation*}
\begin{split}
\int_{\widetilde{M}}\phi ^{2}|\overline{\nabla }_{b}|\overline{\nabla }%
_{b}v||^{2}|\overline{\nabla }_{b}v|^{q-2}d\overline{V}& \leq
-C_{q}^{-1}\int_{\widetilde{M}}\phi ^{2}\overline{\Delta }_{b}|\overline{%
\nabla }_{b}v|^{q}d\overline{V} \\
& =-C_{q}^{-1}\int_{\widetilde{M}}\overline{\nabla }_{b}(\phi ^{2})\overline{%
\nabla }_{b}|\overline{\nabla }_{b}v|^{q}d\overline{V} \\
& \leq 2qC_{q}^{-1}\int_{\widetilde{M}}|\phi ||\overline{\nabla }_{b}\phi ||%
\overline{\nabla }_{b}v|^{q-1}|\overline{\nabla }_{b}|\overline{\nabla }%
_{b}v||d\overline{V} \\
& \leq C_{q}^{\prime }\int_{\widetilde{M}}|\overline{\nabla }_{b}\phi |^{2}|%
\overline{\nabla }_{b}v|^{q}d\overline{V}
\end{split}%
\end{equation*}%
\noindent for some constant $C_{q}^{\prime },$ where the last inequality is
deduced by applying the Schwarz (or Young's) inequality with $\varepsilon $.
Now we change the integral on the right hand side to a corresponding one
using the form $\theta $ and get 
\begin{equation*}
\begin{split}
\int_{\widetilde{M}}\phi ^{2}|\overline{\nabla }_{b}|\overline{\nabla }%
_{b}v||^{2}|\overline{\nabla }_{b}v|^{q-2}d\overline{V}& \leq C_{q}^{\prime
}\int_{\widetilde{M}}\overline{G}^{\frac{-2}{n}}|{\nabla }_{b}\phi |^{2}|{%
\nabla }_{b}v|^{q}\overline{G}^{\frac{-q}{n}}(\overline{G}^{\frac{2}{n}%
})^{n+1}dV \\
& =C_{q}^{\prime }\int_{\widetilde{M}}|\nabla _{b}\phi |^{2}|\nabla
_{b}v|^{q}\overline{G}^{\frac{2(n+1)-2-q}{n}}dV.
\end{split}%
\end{equation*}

\noindent Observe that $q=\frac{2n}{n+1}$ implies $\frac{2(n+1)-2-q}{n}=q$.
We have shown the desired inequality (\ref{8}) for $\phi \in C_{0}^{\infty }(%
\widetilde{M}\setminus \Phi ^{-1}(y))$.\newline
Now for $\phi \in C_{0}^{\infty }(\widetilde{M})$, we consider $\psi
_{r}\phi $ where $\psi _{r}$ is a cutoff function such that for each $\bar{p}
$ $\in $ $\Phi ^{-1}(y),$ $\psi _{r}\equiv 0$ in $B_{r}(\bar{p}),$ $\psi
_{r}\equiv 1$ on $\widetilde{M}\setminus B_{2r}(\bar{p})$, and $0\leq \psi
_{r}\leq 1$ ($r$ small so that \noindent $B_{2r}(\bar{p}_{1})$ $\cap $ $%
B_{2r}(\bar{p}_{2})$ is empty for any pair of points $\bar{p}_{1},$ $\bar{p}%
_{2}$ $\in $ $\Phi ^{-1}(y))$. We also require that $|\nabla _{b}\psi
_{r}|\leq 2r^{-1}$. Applying (\ref{8}) for $\psi _{r}\phi $, we have

\begin{equation*}
\begin{split}
& \int_{\widetilde{M}}\psi _{r}^{2}\phi ^{2}|\overline{\nabla }_{b}|%
\overline{\nabla }_{b}v||^{2}|\overline{\nabla }_{b}v|^{q-2}d\overline{V} \\
& \leq C_{1}\int_{\widetilde{M}}\psi _{r}^{2}|\nabla _{b}\phi |^{2}|\nabla
_{b}v|^{q}\overline{G}^{q}dV+C_{2}\int_{\widetilde{M}}\phi ^{2}|\nabla
_{b}\psi _{r}|^{2}|\nabla _{b}v|^{q}\overline{G}^{q}dV.
\end{split}%
\end{equation*}%
\noindent Noticing that the last integral has order $O(r^{2n-q})\rightarrow
0 $, as $r\rightarrow 0$, we see that (\ref{8}) holds for $\phi \in
C_{0}^{\infty }(\widetilde{M})$.\newline
Next we are going to prove 
\begin{equation}
\begin{split}
& \int_{B_{\rho }(p)}|\overline{\nabla }_{b}v|^{q-2}|\overline{\nabla }_{b}|%
\overline{\nabla }_{b}v|^{2}d\overline{V} \\
& \leq C\rho ^{-2}\int_{B_{4\rho }(p)\backslash B_{\rho
/2}(p)}G^{q}(1+|\nabla _{b}\log {G}|^{2})dV,
\end{split}
\label{19}
\end{equation}%
\noindent where $\rho >0$ is sufficiently large and $C$ is a constant. We
want to make use of (\ref{8}). Note that

\begin{equation}
\begin{split}
\overline{G}^{q}|\nabla _{b}v|^{q}& =|\overline{G}\nabla _{b}v|^{q} \\
& =|\nabla _{b}G-G\overline{G}^{-1}\nabla _{b}\overline{G}|^{q} \\
& \leq C(|\nabla _{b}G|^{q}+G^{q}|\nabla _{b}\log {\overline{G}}|^{q}).
\end{split}
\label{10}
\end{equation}%
\noindent Thus, if we take $\phi \in C_{0}^{\infty }(\widetilde{M})$ such
that $\phi \equiv 1$ in $B_{\rho }(p),\phi \equiv 0$ on $\widetilde{M}%
\setminus B_{2\rho }(p),0\leq \phi \leq 1$ and $|\nabla _{b}\phi |\leq 2\rho
^{-1}$, We see from (\ref{8}) and (\ref{10}) that

\begin{equation}
\begin{split}
\int_{B_{\rho }(p)}|\overline{\nabla }_{b}v|^{q-2}|\overline{\nabla }_{b}|%
\overline{\nabla }_{b}v||^{2}d\overline{V}& \leq \int_{\widetilde{M}}\phi
^{2}|\overline{\nabla }_{b}v|^{q-2}|\overline{\nabla }_{b}|\overline{\nabla }%
_{b}v||^{2}d\overline{V} \\
& \leq C\rho ^{-2}\int_{B_{2\rho }(p)\backslash B_{\rho }(p)}(|\nabla
_{b}G|^{q}+G^{q}|\nabla _{b}\log {\overline{G}}|^{q})dV.
\end{split}
\label{11}
\end{equation}

Let $a$ $=$ $\frac{(2-q)q}{2}.$ Note that $q$ $<$ $2.$ By Young's inequality
we have

\begin{equation}
\begin{split}
|\nabla _{b}G|^{q}& =G^{a}G^{-a}|\nabla _{b}G|^{q} \\
& \leq C(G^{a\frac{2}{2-q}}+G^{-a\frac{2}{q}}|\nabla _{b}G|^{2}) \\
& \leq C(G^{q}+G^{q-2}|\nabla _{b}G|^{2}) \\
& =CG^{q}(1+|\nabla _{b}\log {G}|^{2}),
\end{split}
\label{12}
\end{equation}%
\noindent and 
\begin{equation}
\begin{split}
G^{q}|\nabla _{b}\log {\overline{G}}|^{q}& =G^{q}\frac{|\nabla _{b}\overline{%
G}|^{q}}{\overline{G}^{q}} \\
& \leq G^{q}\frac{C\overline{G}^{q}(1+|\nabla _{b}\log {\overline{G}}|^{2})}{%
\overline{G}^{q}} \\
& =CG^{q}(1+|\nabla _{b}\log {\overline{G}}|^{2}).
\end{split}
\label{13}
\end{equation}

\noindent So from (\ref{11}), (\ref{12}), and (\ref{13}), we obtain 
\begin{equation}
\begin{split}
& \int_{B_{\rho }(p)}|\overline{\nabla }_{b}v|^{q-2}|\overline{\nabla }_{b}|%
\overline{\nabla }_{b}v||^{2}d\overline{V} \\
\leq & C\rho ^{-2}\int_{B_{2\rho }\backslash B_{\rho }}G^{q}(1+|\nabla
_{b}\log {G}|^{2}+|\nabla _{b}\log {\overline{G}}|^{2})dV.
\end{split}
\label{14}
\end{equation}

\noindent Taking $\phi $ $=$ $\psi G^{\frac{q}{2}}$ in Lemma 4.2, where $%
\psi \in C_{0}^{\infty }(\widetilde{M}\setminus \{p\})$, we have 
\begin{equation}
\begin{split}
& \int_{\widetilde{M}}\psi ^{2}G^{q}|\nabla _{b}\log {\overline{G}}|^{2}dV \\
\leq & C\int_{\widetilde{M}}(\psi ^{2}G^{q}|\nabla _{b}\log {G}|^{2}+\nabla
_{b}(\psi G^{\frac{q}{2}})|^{2})dV \\
\leq & C\int_{\widetilde{M}}(\psi ^{2}|\nabla _{b}\log {G}|^{2}+|\nabla
_{b}\psi |^{2})G^{q}dV.
\end{split}
\label{18}
\end{equation}

\noindent Choosing the cutoff function $\psi $ appropriately in (\ref{18}),
we get%
\begin{eqnarray}
&&\int_{B_{2\rho }\backslash B_{\rho }}G^{q}|\nabla _{b}\log {\overline{G}}%
|^{2}dV  \label{18.1} \\
&\leq &C\int_{B_{4\rho }(p)\backslash B_{\rho /2}(p)}G^{q}(1+|\nabla
_{b}\log {G}|^{2})dV  \notag
\end{eqnarray}

\noindent for $\rho $ large. Substituting (\ref{18.1}) into (\ref{14}) gives
(\ref{19}).\newline

Since $b_{n}\Delta _{b}G+R_{\theta }G=0$ on $\widetilde{M}\setminus \{p\}$
and $R_{\theta }\geq C>0$, we have%
\begin{equation*}
\Delta _{b}G=-b_{n}^{-1}R_{\theta }G\leq -b_{n}^{-1}CG
\end{equation*}%
\noindent Multiplying by $\phi ^{2}G^{q-1}$ with $\phi \in C_{0}^{\infty }(%
\widetilde{M}\setminus \{p\})$ and integrating by parts give 
\begin{equation}
\begin{split}
0& \geq \int_{\widetilde{M}}\phi ^{2}G^{q-1}(\Delta _{b}G)dV \\
& =\int_{\widetilde{M}}\nabla _{b}(\phi ^{2}G^{q-1})\nabla _{b}GdV \\
& =\int_{\widetilde{M}}2\phi G^{q-1}\nabla _{b}\phi \nabla
_{b}GdV+(q-1)\int_{\widetilde{M}}\phi ^{2}G^{q-2}|\nabla _{b}G|^{2}dV.
\end{split}
\label{21}
\end{equation}%
\noindent Applying the Schwarz inequality with $\varepsilon $ to (\ref{21}),
we obtain 
\begin{equation}
\int_{\widetilde{M}}\phi ^{2}G^{q-2}|\nabla _{b}G|^{2}dV\leq C\int_{%
\widetilde{M}}G^{q}|\nabla _{b}\phi |^{2}dV.  \label{21.1}
\end{equation}%
\noindent Noting that $G^{q-2}|\nabla _{b}G|^{2}$ $=$ $G^{q}|\nabla _{b}\log
G|^{2}$ and choosing some appropriate $\phi $ in (\ref{21.1}), we can reduce
(\ref{19}) to 
\begin{equation}
\int_{B_{\rho }(p)}|\overline{\nabla }_{b}v|^{q-2}|\overline{\nabla }_{b}|%
\overline{\nabla }_{b}v||^{2}dV\leq C\rho ^{-2}\int_{\widetilde{M}\backslash
B_{\rho /4}}G^{q}dV.  \label{22}
\end{equation}

\noindent By (\ref{3.12.1}) $G^{q}$ is integrable since $q$ $=$ $\frac{2n}{%
n+1}$ $>$ $1$ for $n$ $\geq $ $2.$ So letting $\rho \rightarrow \infty $ in (%
\ref{22}) we get%
\begin{equation*}
|\overline{\nabla }_{b}v|=const.
\end{equation*}

\noindent Since \TEXTsymbol{\vert}$\overline{\nabla }_{b}v|$ $=$ $\bar{G}%
^{-1}|\overline{\nabla }_{b}G-v\overline{\nabla }_{b}\bar{G}|$ $\rightarrow $
$0$ at $p$ by Lemma 4.1, we have $\overline{\nabla }_{b}v$ $=$ $0.$ So $v$ $%
= $ $const..$ From $v\rightarrow 1$ at $p,$ we conclude that $v$ $\equiv $ $%
1.$ \endproof

\section{The positive CR mass theorem}

\label{seccrmass} In this section, according to the work of Li in \cite{Li3}%
, we would like to introduce a positive mass theroem for spherical CR
manifolds. Let $M$ be a closed spherical CR manifold, $\widetilde{M}$ be its
universal cover and 
\begin{equation*}
\pi : \widetilde{M}\longrightarrow M
\end{equation*}
be the canonical projection map, 
\begin{equation*}
\Phi : \widetilde{M}\longrightarrow S^{2n+1}
\end{equation*}
be a CR developing map. We would like to construct local coordinates near
each point $b$ of $M$. There is a local inverse $\pi^{-1}$ as follows: 
\begin{equation*}
\pi^{-1} : U_{b} \longrightarrow \widetilde{M},
\end{equation*}
where $U_{b}$ is a neighborhood of $b\in M$. Let $q=\Phi(p)\in S^{2n+1}$,
where $p\in\pi^{-1}(b)$, the local CR transformation 
\begin{equation*}
T=C_{q}\circ \Phi\circ\pi^{-1} : U_{b} \longrightarrow H^{n}
\end{equation*}
provides $M$ a local coordinate $(z,t)\in H^{n}$ such that $%
(z(b),t(b))=\infty$. Here $C_{q} : S^{2n+1}\longrightarrow H^{n}$ is the
Cayley transform with pole at $q$, i.e. $C_{q}(q)=\infty$. We will call such
coordinates "spherical CR coordinates near $\infty$".

Let $G_{b}$ be the Green's function of $D_{\theta}$ with pole at $b$. It
follows that there is a positive smooth function $h=h(z,t)$ defined on $%
H^{n} $ near $\infty$ such that 
\begin{equation}  \label{secequ8}
(T^{-1})^{*}(G_{b}^{\frac{2}{n}}\theta)=h^{\frac{2}{n}}\Theta.
\end{equation}
By positive constant rescaling we may assume that the complex Jacobian at $p$
is $|\Phi^{^{\prime }}(p)|=1$. Let $\rho(z,t)=(|z|^{4}+t^{2})^{1/4}$ be the
Heisenberg norm on $H^{n}$. Therefore, We have the following asymptotic
expansion of $h=h(z,t)$ near $\infty$:

\bigskip

\textbf{Lemma 5.1.} \textit{Let} $M$ \textit{\ be a closed spherical CR
manifold which is not the standard sphere. Suppose the CR developing map is
injective. Let} $h$ \textit{\ be defined as above. Then we have, near} $%
\infty$, 
\begin{equation*}
h=h(z,t)>1
\end{equation*}
\textit{and } 
\begin{equation*}
h(z,t)=1+A_{b}\cdot\rho(z,t)^{-2n}+O(\rho(z,t)^{-2n-1}).
\end{equation*}

\bigskip

\begin{proof}
Since the projection $\pi$ doesn't change geometry, it follows that near $%
x\in \pi^{-1}(b)$: 
\begin{equation*}
D_{\theta}(\pi^{*}G_{b})=D_{\theta}(G_{b}),
\end{equation*}
where the left hand side is over $\widetilde{M}$ and the right hand side is
over $M$. For $b\in M$ and $x\in\pi^{-1}(b)$, let $\delta_{x}$ be the Dirac
delta function with pole at $x$. Therefore, 
\begin{equation*}
D_{\theta}(\pi^{*}G_{b})=\sum_{x\in\pi^{-1}(b)}\delta_{x},
\end{equation*}
so $\pi^{*}G_{b}$ has poles precisely in the set $\pi^{-1}(b)\subset 
\widetilde{M}$. For each fixed $p\in\pi^{-1}(b)$, 
\begin{equation*}
D_{\theta}(\overline{G}_{p})=\delta_{p}.
\end{equation*}
Therefore $\pi^{*}G_{b}-\overline{G}_{p}$ is bounded near $p$. On the other
hand, by the normalization $|\Phi^{^{\prime }}(p)|=1$ and $\Phi $ is
injective, we have 
\begin{equation*}
\overline{G}_{p}^{\frac{2}{n}}\theta=\Phi^{*}(H_{q}^{\frac{2}{n}%
}\theta_{S^{2n+1}}),
\end{equation*}
it follows that globally on $\widetilde{M}\setminus p$ 
\begin{equation}  \label{secequ12}
((C_{q}\circ\Phi)^{-1})^{*}(\overline{G}_{p}^{\frac{2}{n}}\theta)=\Theta.
\end{equation}
From equations (\ref{secequ8}) and (\ref{secequ12}) and that $\pi^{*}G_{b}$
and $\overline{G}_{p}$ have the same pole strength near $p$, 
\begin{equation*}
h(\infty)=\lim_{(z,t)\rightarrow \infty}h(z,t)=1.
\end{equation*}
On the other hand, because $\overline{G}_{p}$ is the minimal Green's
function for $D_{\theta}$ on $\widetilde{M}$, Bony's strong maximum
principle implies that globally on $\widetilde{M}$: 
\begin{equation*}
\pi^{*}G_{b}>\overline{G}_{p},
\end{equation*}
which by equations (\ref{secequ8}) and (\ref{secequ12}) that for $(z,t)$
near $\infty$, we have $h=h(z,t)>1$.

Next, we would like to show the asymptotic expansion of the function $%
h=h(z,t)$ near $\infty$. By the transformation rule of the CR invariant
sublaplacian, we have $R(G_{b}^{\frac{2}{n}}\theta)=0$ globally on $%
M\setminus \{b\}$, where $R(G_{b}^{\frac{2}{n}}\theta)$ is the Webster
curvature with respect to the contact form $G_{b}^{\frac{2}{n}}\theta$.
Using the transformation rule again we have $D_{\Theta}(h)=0$ for $(z,t)$ in
a neighborhhod of $\infty\in H^{n}$. It follows that (since $h$ is regular
near $\infty$ and $h(\infty)=1$): 
\begin{equation*}
h(z,t)=1+c_{-1}\cdot\rho^{-1}+\cdots +A_{b}\cdot\rho^{-2n}+O(\rho^{-2n-1}).
\end{equation*}
Consider the global CR inversion $\vartheta : H^{n}\setminus
0\longrightarrow H^{n}\setminus 0$ as follows: 
\begin{equation*}
\vartheta(z,t)=(\hat{z},\hat{t})=(-z/w, t/|w|^{2}),
\end{equation*}
where $w=t+i|z|^{2}$. Consider the following standard contact form 
\begin{equation*}
\Theta(z,t)=dt-i\sum_{\alpha=1}^{n}(z^{\alpha}\cdot
dz^{\bar\alpha}-z^{\bar\alpha}\cdot dz^{\alpha}),
\end{equation*}
and 
\begin{equation*}
\Theta(\hat{z},\hat{t})=d\hat{t}-i\sum_{\alpha=1}^{n}(\hat{z}^{\alpha}\cdot d%
\hat{z}^{\bar\alpha}-\hat{z}^{\bar\alpha}\cdot d\hat{z}^{\alpha}).
\end{equation*}
It follows that $\rho(z,t)^{2}=|w|=|\hat{w}|^{-1}=\rho(\hat{z},\hat{t})^{-2}$
and 
\begin{equation*}
\Theta(\hat{z},\hat{t})=|w|^{-2}\cdot \Theta(z,t)=(\rho(\hat{z},\hat{t}%
)^{2n})^{\frac{2}{n}}\cdot\Theta(z,t).
\end{equation*}
Therefore, 
\begin{equation*}
h(z,t)^{\frac{2}{n}}\cdot \Theta(z,t)=(h(\hat{z},\hat{t})\cdot\rho(\hat{z},%
\hat{t})^{-2n})^{\frac{2}{n}}\cdot\Theta(\hat{z},\hat{t}).
\end{equation*}
By the facts that $R(h(z,t)\cdot\Theta(z,t))=0$ and $\Theta(\hat{z},\hat{t})$
and the trasformation rule of the CR invariant sublaplacian, it follows that
near but away from the origin in $H^{n}$, we have 
\begin{equation*}
D_{\Theta(\hat{z},\hat{t})}(h(\hat{z},\hat{t})\cdot\rho(\hat{z},\hat{t}%
)^{-2n})=0, \text{and}\ D_{\Theta(\hat{z},\hat{t})}\rho(\hat{z},\hat{t}%
)^{-2n}=0.
\end{equation*}
Therefore, near but away fron the origin in $H^{n}$: 
\begin{equation*}
D_{\Theta(\hat{z},\hat{t})}(h(\hat{z},\hat{t})\cdot\rho(\hat{z},\hat{t}%
)^{-2n}-\rho(\hat{z},\hat{t})^{-2n})=0
\end{equation*}
and note that 
\begin{equation*}
\begin{split}
&h(\hat{z},\hat{t})\cdot\rho(\hat{z},\hat{t})^{-2n}-\rho(\hat{z},\hat{t}%
)^{-2n} \\
=&c_{-1}\cdot\rho(\hat{z},\hat{t})^{-2n+1}+\cdots +c_{-2n+1}\cdot\rho(\hat{z}%
,\hat{t})^{-1}+A_{b}+O(\rho(\hat{z},\hat{t})).
\end{split}%
\end{equation*}
A standard removable singularity argument (Proposition 5.17 in \cite{JL})
implies that 
\begin{equation*}
c_{-1}=\cdots =c_{-2n+1}=0.
\end{equation*}
Therefore, we have the following aymptotic expansion of $h=h(z,t)$ near $%
\infty$: 
\begin{equation*}
h(z,t)=1+A_{b}\cdot\rho(z,t)^{-2n}+O(\rho(z,t)^{-2n-1}).
\end{equation*}
This completes the lemma.
\end{proof}

\begin{de}
We call the constant $A_{b}$ CR mass.
\end{de}

We would like to remark that the constant $A_b$ doesn't depend on the choice
of local coordinates near $b\in M$ (see \cite{Li3}). Corollary C states the
positivity of the CR mass.

\proof
\textbf{(of Corollary C)}

By Theorem A, $\Phi$ is injective. It follows from Lemma 5.1 that for $(z,t)$
near $\infty$ 
\begin{equation*}
h=h(z,t)=1+A_{b}\cdot\rho(z,t)^{-2n}+O(\rho(z,t)^{-2n-1})>1.
\end{equation*}
Let $B_{L}(0)$ be a ball on $H^{n}$ centered at $0$ and with radius $L$ such
that 
\begin{equation*}
m=\min_{\partial B_{L}(0)}(h-1)>0.
\end{equation*}
Since in $H^{n}\setminus B_{L}(0)$, 
\begin{equation*}
D_{\Theta}(h-1)=D_{\Theta}(mL^{2n}\rho^{-2n})=0
\end{equation*}
and 
\begin{equation*}
(h-1)|_{\partial B_{L}(0)}\geq mL^{2n}\rho^{-2n}|_{\partial B_{L}(0)},
\end{equation*}
we conclude from Bony's maximum principle that in $H^{n}\setminus B_{L}(0),$ 
\begin{equation*}
h-1\geq mL^{2n}\rho^{-2n}.
\end{equation*}
Therefore, $A_{b}\geq mL^{2n}>0$. \endproof

\endproof%


\bigskip

\begin{thebibliography}{99}
\bibitem{A} Aubin, T., Nonlinear Analysis on Manifolds. Monge-Amp\'{e}re
Equations. \textit{Springer Verlag, }Berlin and New York\textit{, }1982.

\bibitem{Bony} Bony, J. M., Principe du maximu, inegalite de Harnack, et
unicite du probleme de Cauchy pour le operateurs elliptiques de degeneres, 
\textit{Ann. Inst. Fourier (Grenoble) }\textbf{19} (1) (1969), 277-304.

\bibitem{BS} Burns, D. and Shnider, S., Spherical hypersurfaces in complex
manifolds, \textit{Invent. Math., }\textbf{33} (1976), 223-246\textit{.}

\bibitem{Chiu} Chiu, H.-L., The sharp lower bound for the first positive
eigenvalue of the sublaplacian on a pseudohermitian 3-manifold, \textit{Ann.
Glob. Anal. Geom.}, \textbf{30} (2006), 81-96.

\bibitem{CC} Chang, S.-C. and Chiu, H.-L., Nonnegativity of $CR$ Paneitz
operator and its application to the $CR$ Obata's theorem, \textit{J. Geom.
Anal}., \textbf{19} (2009), 261-287.

\bibitem{CCY1} Chanillo, S., Chiu, H.-L., and Yang, P., Embeddability for
three-dimensional Cauchy-Riemann manifolds and CR Yamabe invariants, \textit{%
to appear in Duke Math. J.}.

\bibitem{CCY2} Chanillo, S., Chiu, H.-L., and Yang, P., Nonnegativity
criterion for Paneitz operators on three-dimensional $CR$ manifolds, \textit{%
preprint.}

\bibitem{CMY} Cheng, J.-H., Malchiodi, A., and Yang, P., A positive mass
theorem in Cauchy-Riemann geometry of dimension 3, \textit{preprint}.

\bibitem{CT} Cheng, J.-H. and Tsai, I-H., Deformation of spherical $CR$
structures and the universal Picard variety, \textit{Commun. in Anal. and
Geom.}, \textbf{8} (2000), 301-346.

\bibitem{CM} Chern, S.-S. and Moser, J., Real hypersurfaces in complex
manifolds, \textit{Acta Math}., \textbf{133} (1974), 219-271.

\bibitem{FG} Falbel, E. and Gusevskii, N., Spherical $CR$-manifolds of
dimension 3, \textit{Bol. Soc. Brasil. Mat. (N.S.)} \textbf{25} (1994),
31-56.

\bibitem{FS} Folland, G. B. and Stein, E. M., Estimates for the $\bar{%
\partial}_{b}$-complex and analysis on the Heisenberg group, \textit{Comm.
Pure Appl. Math.}, \textbf{27} (1974), 429-522.

\bibitem{GL} Graham, C. R. and Lee, J. M., Smooth solutions of degenerate
Laplacians on strictly pseudoconvex domains, \textit{Duke Math. J.}, \textbf{%
57} (1988), 697-720.

\bibitem{Gol} Goldman, W. M., Complex hyperbolic geometry, \textit{Oxford
Mathematical Monographs, Clarendon Press/Oxford University Press}, New York,
1999.

\bibitem{Gr} Greenleaf, A., The first eigenvalue of a sublaplacian on a
pseudohermitian manifold, \textit{Commun. Partial Differential Equations} 
\textbf{10(2)} (1985), 191-217.

\bibitem{Hir} Hirachi, K., Scalar pseudo-Hermitian invariants and the Szeg%
\"{o} kernel on 3-dimensional $CR$ manifolds, \textit{Complex geometry,
Lecture Notes in Pure and Appl. Math.}, Vol. \textbf{143}, pp 67-76, Dekker,
New York, 1992.

\bibitem{JL} Jerison, D. and Lee, J. M., The Yamabe problem on $CR$
manifolds, \textit{J. Diff. Geom.} \textbf{25} (1987), 167-197.

\bibitem{KT} Kamishima, Y. and Tsuboi, T., $CR$-structures on Seifert
manifolds, \textit{Invent. Math}. \textbf{104} (1991), 149-163.

\bibitem{Lee} Lee, John M., The Fefferman metric and pseudohermitian
invariants, \textit{Trans. Amer. Math. Soc}., \textbf{296} (1986), 411-429.

\bibitem{Lee2} Lee, John M., Pseudo-Einstein structures on $CR$ manifolds, 
\textit{Amer. J. Math.}, \textbf{110} (1988), 157-178.

\bibitem{JL} Jerison, D. and Lee, J. M., The Yamabe problem on $CR$
manifolds, \textit{J. Diff. Geom.}, \textbf{25} (1987), 167-197.

\bibitem{Li1} Li, Z., On spherical $CR$ manifolds with positive Webster
scalar curvature. Preprint.

\bibitem{Li3} Li, Z., The Yamabe problem on spherical $CR$ manifolds of
dimensions $\geq $ 7, Preprint.

\bibitem{Sch} Schwartz, R. E., Spherical $CR$ geometry and Dehn surgery, 
\textit{Princeton University Press, }Princeton and Oxford, 2007

\bibitem{Str} Strichartz, R. S., Sub-Riemannian geometry, \textit{J. Diff.
Geom.}, \textbf{24} (1986), 221-263.

\bibitem{SY1} Schoen, R. and Yau, S.-T., Conformally flat manifolds,
Kleinian groups and scalar curvature. \textit{Invent. Math.,} \textbf{92}
(1988), 47-71.

\bibitem{SY2} Schoen, R. and Yau, S.-T., Lectures on Differential Geometry,
Conference Proceedings and Lecture Notes in Geometry and Topology, Vol 1, 
\textit{International Press}, 1994.

\bibitem{Ta} Tanaka, N., A differential geometric study on strongly
pseudoconvex manifolds, \textit{Kinokuniya Company Ltd., Tokyo}, 1975.

\bibitem{We} Webster, S., Pseudohermitian structures on a real hypersurface, 
\textit{J. Diff. Geom}. \textbf{13} (1978), 25-41.
\end{thebibliography}
\end{document}